\documentclass[]{article}
\begin{document}
\newtheorem{proposition}{Proposition}[section]
\newtheorem{definition}{Definition}[section]
\newtheorem{lemma}{Lemma}[section]

\title{\bf Automatic Integration  }
\author{Keqin Liu\footnote{ Department of Mathematics, The University of British Columbia, Vancouver, BC,
Canada, V6T 1Z2,  Email address: kliu at math.ubc.ca}}
\date{June, 2020}
\maketitle

\begin{abstract}  The purpose of this paper is to introduce  the concept of the  automatic integration and present a new way of  approximating  definite integrals using   the  automatic integration based on an associative algebra with zero divisors.

\end{abstract}

\medskip
Both numerical  differentiation and symbolic differentiation have been replaced by automatic differentiation  extensively in scientific computation, especially, in machine learning community. The motivation of doing 
automatic differentiation comes from the unwonted fact: a satisfying method of evaluating derivatives to meet the need of large-scale machine learning is not found even many different efforts were made  for several decades by developing various numerical methods based on the definition of derivatives and by  finding better
ways of performing  symbolic  differentiation. The key idea which makes automatic differentiation so successful is that both the  numerical  approach of using the  definition of derivatives and the symbolic approach of using explicit formulas of derivatives should be replaced by the algebraic 
approach which is based on  the  strategy of  computing the value of a function from an  associative algebra with zero-divisors to the algebra itself. 

\medskip
Evaluating  derivatives and evaluating definite integrals are two fundamental problems. 
In many applications from engineering and statistics, 
the main tool of evaluating one-dimensional definite integrals is numerical integration or numerical quadrature. During the long history of numerical integration, various techniques 
including midpoint rule, trapezoid rule and Simpson's rule have been developed to do  numerical integration. However, all of these techniques are not at all connected to the key idea which makes automatic differentiation to gain the great success in scientific computation. A freshman knows that an anti-derivative of a function is enough to get the definite integral of the function.  Also, it is well-known that finding an anti-derivative of a function is much more difficult than  evaluating  the derivative of the function. The reason for this well-known fact is that you can almost always get the derivative of a differential function as long as you have patient to use the rules of differentiation, but the
anti-derivatives of a great many functions appearing in engineering and statistics can not be obtained usually  no matter how hard you try to use various integration techniques. Based on the fact that the satisfying  solution to the problem of evaluating  derivatives has to use automatic differentiation, it is reasonable to believe that 
the satisfying  solution to the harder problem of evaluating definite integrals  must depend on the new strategy of utilizing the key idea from automatic differentiation.  Replacing the integrand of a definite integral with its higher order Taylor polynomial is a very simple and natural  idea of approximating the definite integral, but this simple and natural  idea has not been used successfully to develop an effective approximation method in numerical integration.
This is clearly due to the fact that calculating  higher order derivatives is both quite complicated and over-elaborate if we just use  numerical  differentiation and symbolic differentiation.  With the advent of the success of automatic differentiation in scientific computation,  it is time to combine  the simple and natural  idea of using Taylor polynomials with automatic differentiation technique to see if  an effective method of approximating definite integrals can be obtained.  The purpose of this paper is to initiate the study of this new approach.  To introduce this new approach clearly, we will explain how to combine  the  idea of using $5$th-order Taylor polynomials with automatic differentiation technique to approximate definite integrals in this paper.

\medskip
This paper consists of three sections. In section 1,  we define ${k_1, k_2, \dots , k_N \choose  x_1, x_2, \cdots , x_{N-1}}$-automatic approximation of  the definite integral $\displaystyle\int_a^b f(x) dx$ at its  center 
\linebreak  $(c_1, c_2, \dots , c_N)$,  introduce  the concept of the $n$th-order automatic integration
and give the $\mathcal{R}^{(6)}$-extensions of some common elementary functions, where 
$\mathcal{R}^{(6)}$, which is denoted by $\mathcal{P}_{6}$ in  Section 13.2 of \cite{GW}, is the  $6$-dimensional truncated polynomial real algebra  $\mathcal{R}^{(6)}$. In section 2, we prove the main theorem of this paper which presents  the  $5$th-order  automatic integration technique of  computing $5$-automatic approximation of  the definite integral $\displaystyle\int_a^b f(x) dx$ at its  center  $c$.  In section 3, as an example, we compute the  ${k_1, k_2, \dots , k_N \choose  x_1, x_2, \cdots , x_{N-1}}$-automatic approximations of  the definite integral $\displaystyle\int_0^2 e^{x^2} dx$ at different  centers   $(c_1, c_2, \dots , c_N)$ for $k_1= k_2=\cdots = k_N=5$ and $1\le N\le 3$. There are two obvious facts  which appear in our computation of the  ${k_1, k_2, \dots , k_N \choose  x_1, x_2, \cdots , x_{N-1}}$-automatic approximations of  the definite integral $\displaystyle\int_0^2 e^{x^2} dx$ at different  centers   $(c_1, c_2, \dots , c_N)$. One  fact is that although the accuracy of  the Midpoint Rule,  the Trapezoid Rule and Simpson approximations to the definite integral $\displaystyle\int_0^2 e^{x^2} dx$ just depend on the number of the subintervals of $[0, \,2]$,   the accuracy of  the  ${k_1, k_2, \dots , k_N \choose  x_1, x_2, \cdots , x_{N-1}}$-automatic approximations of  the definite integral $\displaystyle\int_0^2 e^{x^2} dx$ at different  centers   $(c_1, c_2, \dots , c_N)$  not only depend on the the number of the subintervals of $[0, \,2]$, but also depend on the choices of the centers   $(c_1, c_2, \dots , c_N)$. The other fact is that to get the same accuracy, the number of  the subintervals of $[0, \,2]$ used in the  ${k_1, k_2, \dots , k_N \choose  x_1, x_2, \cdots , x_{N-1}}$-automatic approximations of  the definite integral $\displaystyle\int_0^2 e^{x^2} dx$ can be much less than the number of  the subintervals of $[0, \,2]$ used in the Midpoint Rule,  the Trapezoid Rule and Simpson approximations to the definite integral $\displaystyle\int_0^2 e^{x^2} dx$.

\bigskip
Throughout this paper,   the real number field is denoted by  $\mathcal{R}$ , the range of a function $f$ is denoted by $Im f$, and an associative algebra with the identity is just called  a unital  associative algebra.

\bigskip
\section{ The Concept of Automatic Integration }

Let $[a, \, b]$ be a  closed real number interval, and  let  ${\bf D}^{\infty}[a, \, b]$ be the associative algebra of all analytic functions on an open   interval $I$  which  contains the closed real number interval $[a, \, b]$. In other words, 
$f\in {\bf D}^{\infty}[a, \, b]$ if and only if the Taylor series of $f$ with its center $c\in I$ exists for all $c\in I$.

\medskip
For  $f(x)\in {\bf D}^{\infty}[a, \, b]$ and  $c\in [a, b]$, the polynomial
\begin{equation}\label{eq1}
T_{n; c}(x):=f(c)+f'(c)(x-c)+\displaystyle\frac{f^{(2)}(c)}{2!}(x-c)^{2}+\cdots +
\displaystyle\frac{f^{(n)}(c)}{n!}(x-c)^{n}
\end{equation}
is called the {\bf $n$th-order Taylor polynomial} for $f(x)$ with its {\bf center} at $c$.
It is well-known that $\displaystyle\int_a^b T_{n; c}(x) dx$ is a good approximation of the definite integral 
$\displaystyle\int_a^b f(x) dx$ if $n$ is large enough. We call  $\displaystyle\int_a^b T_{n; c}(x) dx$ the 
{\bf $n$-automatic approximation} to  $\displaystyle\int_a^b f(x) dx$ at its {\bf center}  $c$.

\medskip
In general, suppose $[x_0, x_1]$, $[x_1, x_2]$, \dots , $[x_{N-1}, x_N]$ are subintervals of $[a, b]$ with
$$
a=x_0< x_1< x_2<\cdots <x_{N-1}< x_N=b,
$$
where $N$ is a positive integer. Let 
$f\in {\bf D}^{\infty}[a, \, b]$.  If 
$c_i\in [x_{i-1}, \, x_i]$ for $1\le i\le N$ and $k_1$, $k_2$, $\cdots$, $k_N$ are positive integers, then
 $$\displaystyle\sum_{i=1}^N\displaystyle\int_{x_{i-1}}^{x_i} T_{k_i; c_i}(x) dx$$ is called 
{\bf ${k_1, k_2, \dots , k_N \choose  x_1, x_2, \cdots , x_{N-1}}$-automatic approximation} of  $\displaystyle\int_a^b f(x) dx$ at its {\bf center}  $(c_1, c_2, \dots , c_N)$, where 
 ${k_1, k_2, \dots , k_N \choose  x_1, x_2, \cdots , x_{N-1}}:=k_1$ if $N=1$.
Clearly, ${k_1, k_2, \dots , k_N \choose  x_1, x_2, \cdots , x_{N-1}}$-automatic approximation 
$\displaystyle\sum_{i=1}^N\displaystyle\int_{x_{i-1}}^{x_i} T_{k_i; c_i}(x) dx$ to  $\displaystyle\int_a^b f(x) dx$ with its  center at 
$(c_1, c_2, \dots , c_N)$ approaches the definite integral  $\displaystyle\int_a^b f(x) dx$ more quickly if 
\linebreak
$\min\{k_1, k_2, \dots , k_N\}$ is large enough and $\max\{x_1-x_0, x_2-x_1, \dots , x_{N}-x_{N-1}\}$ is small enough.

\medskip
Instead of using the anti-derivative of  $T_{n; c}(x)$, our new strategy of computing the 
 $n$-automatic approximation   $\displaystyle\int_a^b T_{n; c}(x) dx$ to  $\displaystyle\int_a^b f(x) dx$ at its center $c$ is to use the value of a  $\mathcal{R}^{(n+1)}$-valued function   to get the exact value of the  $n$-automatic approximation 
  $\displaystyle\int_a^b T_{n; c}(x) dx$ to  $\displaystyle\int_a^b f(x) dx$ at its  center  $c$, where 
  $\mathcal{R}^{(n+1)}$ is  the  $(n+1)$-dimensional truncated polynomial real algebra, which is  defined by
$$
\mathcal{R}^{(n+1)}=\displaystyle\bigoplus_{i=0}^{n} \mathcal{R}\varepsilon ^i, \quad\mathcal{R}\varepsilon ^0=\mathcal{R}, \quad
\varepsilon ^i\cdot\varepsilon ^j
=\left\{\begin{array}{cl}\varepsilon ^{i+j}&
\mbox{if $i+j< n+1$}\\0&\mbox{if $i+j\ge n+1$}\end{array}\right.
$$
for $0\le i, \, j\le n$ .

\bigskip
Let $A$ be a unital associative real algebra.  For a non-empty set $S$,  we use 
$${\bf F}(S, A):= \{f\,|\,\mbox{$ f: S\to A$ is a function } \}$$
to denote the set of the functions from $S$ to $A$. For $f$, $g\in {\bf F}(S, A)$, $r\in \mathcal{R}$ and $x\in S$, we define
$$
(f+g)(x):=f(x)+g(x),\quad (rf)(x):=r\cdot f(x),\quad (f\cdot g)(x):=f(x)g(x).
$$
Then ${\bf F}(S, A)$ is a  unital associative real algebra with respect to the addition,  the scalar multiplication and the product above. The  identity $1_{{\bf F}(S, \, A)}$ of the algebra  ${\bf F}(S, A)$ is the constant function given by 
$$
1_{{\bf F}(S, \, A)} (x):=1\quad \mbox{for $x\in S$,}
$$
where $1$ is the identity of the  unital associative real algebra $A$.

\medskip
For $c\in \mathcal{R}$, we define
$$
{\bf D}^n(c):=\left\{f\,\left|\begin{array}{l}
\mbox{$f$ is a real-valued function defined on an open interval}\\
\mbox{of real numbers and $f$ has the $n$-th derivative at $c$}
\end{array}\right.\right\}.
$$
Clearly, ${\bf D}^n(c)$ is a unital associative real algebra  and ${\bf D}^{\infty}[a, \, b]$ is a subalgebra of 
${\bf D}^n(c)$  for all positive integer $n$ and all $c\in [a, \,b]$.

\bigskip
We now introduce the concept of  $n$th-order automatic integration in the following definition.

\medskip
\begin{definition}\label{def1.1} Let $A$ be a unital associative real algebra,  let $[a, \,b]$ be a closed real number interval and let $n$ be a positive integer. 
A $3$-tuple  $\big(\Lambda,\, \Omega, \,\Gamma\big)$ 
consisting of a map $\Lambda: \displaystyle\bigcup _{x\in \mathcal{R}}{\bf D}^n(x)\to  \displaystyle\bigcup_{S\subseteq A}{\bf F}(S, \,A)$, a map   $\Omega: \mathcal{R}\to A$  and a  map $\Gamma: A\to \mathcal{R}$ is called the {\bf $n$th-order automatic integration} induced by $A$  if the following four conditions are satisfied:
\begin{description} 
\item[(i)] For each $x\in \mathcal{R}$, there exists a subset $A_x\subseteq A$ such that $\Omega(x)\in A_x$,
$Im (\Lambda |{\bf D}^n(x))\subseteq {\bf F}(A_x, \,A)$ and the map  
$\Lambda |{\bf D}^n(x): {\bf D}^n(x)\to {\bf F}(A_x, \,A)$ is a real linear transformation;
\item[(ii)] $\Lambda$  preserves the product at $\Omega(x)$ with $x\in \mathcal{R}$, which means 
\begin{equation}\label{eq32A}
\Lambda(f\cdot g)\big(\Omega(x)\big)=(\Lambda f)\big(\Omega(x)\big)\cdot (\Lambda g)\big(\Omega(x)\big)\quad
\mbox{for $f$, $g\in {\bf D}^n(x)$ };
\end{equation}
\item[(iii)]   $\Lambda$  preserves the composition  at $\Omega(x)$ with $x\in \mathcal{R}$, which means that if  $g\in {\bf D}^n(x)$ and  $f\in {\bf D}^n\big(g(x)\big)$, then $Im \Lambda (g)\subseteq
A_{g(x)}$ and 
\begin{equation}\label{eq32}
 \Lambda(f\circ g)\big(\Omega(x)\big)=\big(\Lambda(f)\circ \Lambda(g)\big)\big(\Omega(x)\big);
\end{equation}
\item[(vi)]  For $f\in {\bf D}^{\infty}[a, \, b]$ and $c\in [a, \,b]$, we have
\begin{equation}\label{eq4}
\Big(\Gamma\circ (\Lambda(f))\circ\Omega\Big)(c)=\displaystyle\int_a^b T_{n; c}(x) dx,
\end{equation}
where  $T_{n; c}(x)$ is the $n$th-order Taylor polynomial for $f$ with its center at $c\in [a, b]$.
\end{description}
\end{definition} 

\medskip
For the sake of simplicity, in this paper, we use the  $5$th-order automatic integration induced by 
$\mathcal{R}^{(6)}$ to explain how to compute the  $5$-automatic approximation 
$\displaystyle\int_a^b T_{5; c}(x) dx$ to the definite integral  $\displaystyle\int_a^b f(x) dx$ at its  center 
$c\in [a, b]$. As a preparation, 
we finish this section by indicating the way of extending some   common elementary functions in 
${\bf D}^5(x)$ to the functions in ${\bf F}(\mathcal{R}^{(6)}_x, \,\mathcal{R}^{(6)})$ with $x\in \mathcal{R}$, where
$$
\mathcal{R}^{(6)}_x:=\{x+a_1\varepsilon +a_2\varepsilon^2+a_3\varepsilon^3+a_4\varepsilon^4+a_5\varepsilon^5\,|\, a_1, a_2, a_3, a_4, a_5\in \mathcal{R}\}
$$

\medskip
For  $f\in {\bf D}^5(x)$, we define the map 
$\Lambda: {\bf D}^5(x)\to{\bf F}(\mathcal{R}^{(6)}_x, \,\mathcal{R}^{(6)})$ by
$$
\Lambda (f)(x+a_1\varepsilon +a_2\varepsilon^2+a_3\varepsilon^3 +a_4\varepsilon^4+a_5\varepsilon^5):=f(x)+a_1f'(x)\varepsilon+
$$
$$
+\Big(a_2f'(x)+\displaystyle\frac{1}{2!} a_1^2f^{(2)}(x)\Big)\varepsilon^2+\Big(a_3f'(x)+a_1a_2f^{(2)}(x)+\displaystyle\frac{1}{3!} a_1^3f^{(3)}(x)\Big)\varepsilon^3+
$$
$$
+\Big[a_4f'(x)+\Big(a_1a_3+\displaystyle\frac{1}{2}a_2^2 \Big)f^{(2)}(x)+
\displaystyle\frac{1}{2}a_1^2a_2 f^{(3)}(x)+\displaystyle\frac{1}{4!}a_1^4 f^{(4)}(x)\Big]\varepsilon^4+
$$
\begin{eqnarray}\label{eq5}
&&+\Big[a_5f'(x)+(a_1a_4+a_2 a_3)f^{(2)}(x)+\qquad\qquad\qquad\qquad\qquad\qquad\qquad\nonumber\\
&&\quad +
\displaystyle\frac{1}{2}(a_1^2a_3+a_1a_2^2) f^{(3)}(x)+\displaystyle\frac{1}{6}a_1^3a_2 f^{(4)}(x)+
\displaystyle\frac{1}{5!}a_1^5 f^{(5)}(x)\Big]\varepsilon^5,
\end{eqnarray}
where $x$, $a_1$, $a_2$, $a_3$, $a_4$, $a_5\in \mathcal{R}$.  We call $\Lambda (f)$
the {\bf $\mathcal{R}^{(6)}$-extension} of $f$, which is also denoted by $\overline{f}$. The $\mathcal{R}^{(6)}$-extensions of some common elementary functions in ${\bf D}^5(x)$ are given as follows:
\begin{eqnarray*}
&\bullet&\displaystyle\frac{1}{x+a_1\varepsilon +a_2\varepsilon^2+a_3\varepsilon^3+a_4\varepsilon^4+a_5\varepsilon^5}=
\displaystyle\frac{1}{x}-\displaystyle\frac{a_1}{x^2}\varepsilon 
+\left(\displaystyle\frac{a_1^2}{x^3}-\displaystyle\frac{a_2}{x^2}\right)\varepsilon^2+\\
&&\quad +\left(-\displaystyle\frac{a_1^3}{x^4}+\displaystyle\frac{2a_1a_2}{x^3}
-\displaystyle\frac{a_3}{x^2}\right)\varepsilon^3
+\left(\displaystyle\frac{a_1^4}{x^5}-\displaystyle\frac{3a_1^2a_2}{x^4}
+\displaystyle\frac{2a_1a_3+a_2^2}{x^3}-\displaystyle\frac{a_4}{x^2}\right)\varepsilon^4+\\
&&\quad +\left(-\displaystyle\frac{a_1^5}{x^6}+\displaystyle\frac{4a_1^3a_2}{x^5}
-\displaystyle\frac{3(a_1^2a_3+a_1a_2^2)}{x^4}+\displaystyle\frac{2(a_1a_4+a_2a_3)}{x^3}
-\displaystyle\frac{a_5}{x^2}\right)\varepsilon^5\\
&&\quad\quad\mbox{for $0\ne x\in \mathcal{R}$}\\
&&\\
&\bullet&\overline{\exp} (x+a_1\varepsilon +a_2\varepsilon^2+a_3\varepsilon^3
+a_4\varepsilon^4+a_5\varepsilon^5)=e^x +
a_1e^x \varepsilon +\left(\displaystyle\frac{a_1^2}{2}+a_2\right)e^x\varepsilon^2+\\
&&\quad +\left(\displaystyle\frac{a_1^3}{6}+a_1a_2+a_3\right)e^x\varepsilon^3
+\left(\displaystyle\frac{a_1^4}{24}+\displaystyle\frac12 a_1^2a_2+a_1a_3
+\displaystyle\frac12 a_2^2+a_4\right)e^x\varepsilon^4+\\
&&\quad+\Big(\displaystyle\frac{1}{120}a_1^5 +\displaystyle\frac{1}{6}a_1^3a_2+
\displaystyle\frac{a_1^2a_3+a_1a_2^2}{2}+a_1a_4+a_2 a_3+a_5
\Big)e^x\varepsilon^5\\
&&\\
&\bullet&\overline{\sin}\, (x+a_1\varepsilon +a_2\varepsilon^2+a_3\varepsilon^3
+a_4\varepsilon^4+a_5\varepsilon^5)=\sin x +
(a_1\cos x)\varepsilon+\\
&&\quad +\left(-\displaystyle\frac{a_1^2}{2}\sin x+a_2\cos x\right)\varepsilon^2+
\left(-\displaystyle\frac{a_1^3}{6}\cos x-a_1a_2\sin x+a_3\cos x\right)\varepsilon^3+\\
&&\quad +\left(\displaystyle\frac{a_1^4}{24}\sin x-\displaystyle\frac12 a_1^2a_2\cos x-
\Big(a_1a_3+\displaystyle\frac12 a_2^2\Big)\sin x+a_4\cos x\right)\varepsilon^4+\\
&&\quad +\Big[\displaystyle\frac{a_1^5}{120}\cos x+\displaystyle\frac16 a_1^3a_2\sin x-
\displaystyle\frac{a_1^2a_3+a_1 a_2^2}{2}\cos x+\\
&&\quad\quad \qquad-(a_1a_4+a_2a_3)\sin x+a_5\cos x\Big]\varepsilon^5\\
&&\\
&\bullet&\overline{\cos}\, (x+a_1\varepsilon +a_2\varepsilon^2+a_3\varepsilon^3
+a_4\varepsilon^4+a_5\varepsilon^5)=\cos x -
(a_1\sin x)\varepsilon+\\
&&\quad +\left(-\displaystyle\frac{a_1^2}{2}\cos x-a_2\sin x\right)\varepsilon^2+
\left(\displaystyle\frac{a_1^3}{6}\sin x-a_1a_2\cos x-a_3\sin x\right)\varepsilon^3+\\
&&\quad+\left(\displaystyle\frac{a_1^4}{24}\cos x+\displaystyle\frac12 a_1^2a_2\sin x-
\Big(a_1a_3+\displaystyle\frac12 a_2^2\Big)\cos x-a_4\sin x\right)\varepsilon^4
+\\
&&\quad+\Big[-\displaystyle\frac{a_1^5}{120}\sin x+\displaystyle\frac16 a_1^3a_2\cos x+
\displaystyle\frac{a_1^2a_3+a_1 a_2^2}{2}\sin x+\\
&&\quad\quad \qquad-(a_1a_4+a_2a_3)\cos x-a_5\sin x\Big]\varepsilon^5\\
&&\\
&\bullet&\overline{\ln}\,(x+a_1\varepsilon +a_2\varepsilon^2+a_3\varepsilon^3
+a_4\varepsilon^4+a_5\varepsilon^5)=
\ln x+\displaystyle\frac{a_1}{x}\varepsilon 
+\left(-\displaystyle\frac{a_1^2}{2x^2}+\displaystyle\frac{a_2}{x}\right)\varepsilon^2+\\
&&\quad +\left(\displaystyle\frac{a_1^3}{3x^3}-\displaystyle\frac{a_1a_2}{x^2}
+\displaystyle\frac{a_3}{x}\right)\varepsilon^3
+\left(-\displaystyle\frac{a_1^4}{4x^4}+\displaystyle\frac{a_1^2a_2}{x^3}-
\displaystyle\frac{2a_1a_3+a_2^2}{2x^2}+\displaystyle\frac{a_4}{x}\right)\varepsilon^4+\\
&&\quad +\left(\displaystyle\frac{a_1^5}{5x^5}-\displaystyle\frac{a_1^3a_2}{x^4}+
\displaystyle\frac{a_1^2a_3+a_1a_2^2}{x^3}-\displaystyle\frac{a_1a_4+a_2a_3}{x^2}
+\displaystyle\frac{a_5}{x}\right)\varepsilon^5\quad \quad\mbox{for $0< x\in \mathcal{R}$}\\
&&\\
&\bullet&\overline{\arctan}\,(x+a_1\varepsilon +a_2\varepsilon^2+a_3\varepsilon^3+a_4\varepsilon^4+a_5\varepsilon^5)=
\arctan x+\displaystyle\frac{a_1}{1+x^2}\varepsilon +\\
&&\quad 
+\left(-\displaystyle\frac{a_1^2x}{(1+x^2)^2} +\displaystyle\frac{a_2}{1+x^2}\right)\varepsilon^2
 +\left(\displaystyle\frac{a_1^3(3x^2-1)}{3(1+x^2)^3}-\displaystyle\frac{2a_1a_2x}{(1+x^2)^2}
+\displaystyle\frac{a_3}{1+x^2}\right)\varepsilon^3+\\
&&\quad +\left(\displaystyle\frac{a_1^4(x-x^3)}{(1+x^2)^4} +
\displaystyle\frac{a_1^2a_2(3x^2-1)}{(1+x^2)^3}-\displaystyle\frac{(2a_1a_3+a_2^2)x}{(1+x^2)^2}
+\displaystyle\frac{a_4}{1+x^2}\right)\varepsilon^4+\\
&&\quad +\left(\displaystyle\frac{a_1^5(1-10x^2+5x^4)}{5(1+x^2)^5} +
\displaystyle\frac{4a_1^3a_2(x-x^3)}{(1+x^2)^4} +\displaystyle\frac{(a_1^2a_3+a_1a_2^2)(3x^2-1)}{(1+x^2)^2}+\right.\\
&&\quad\qquad\qquad  \left. -\displaystyle\frac{(a_1a_4+a_2a_3)x}{(1+x^2)^2}+
\displaystyle\frac{a_5}{1+x^2}\right)\varepsilon^5
\end{eqnarray*}

\bigskip
\section {Automatic integration induced by $\mathcal{R}^{(6)}$}

The following theorem, which is the main theorem of this paper, presents the 
new technique of using the  $5$th-order automatic integration induced by $\mathcal{R}^{(6)}$  to compute the  $5$-automatic approximation 
$\displaystyle\int_a^b T_{5; c}(x) dx$ to the definite integral  $\displaystyle\int_a^b f(x) dx$ at its  center 
$c\in [a, b]$. 

\medskip
\begin{proposition}\label{pr2.1} ({\bf The Main Theorem})  Let $[a, b]$ be a real number interval, and let 
$\beta_1$,  $\dots$,  $\beta_5$ be real constants with 
$\beta_1\ne 0$. 
If the map $\Lambda$  is defined by (\ref{eq5}) and 
the real numbers $A_1$, $\dots$, $A_5$ are given by
\begin{equation}\label{eq8}
\left\{\begin{array}{l}
A_5:=\displaystyle\frac{(b-c)^6-(a-c)^6}{6\beta_1^5}, \\\\
A_4:=\displaystyle\frac{(b-c)^5-(a-c)^5}{5\beta_1^4}-\displaystyle\frac{4\beta_2A_5}{\beta_1}, \\\\
A_3:=\displaystyle\frac{(b-c)^4-(a-c)^4}{4\beta_1^3}-\displaystyle\frac{3\beta_2A_4}{\beta_1}
-3\left(\displaystyle\frac{\beta_3}{\beta_1}+\displaystyle\frac{\beta_2^2}{\beta_1^2}\right)A_5, \\\\
A_2:=\displaystyle\frac{(b-c)^3-(a-c)^3}{3\beta_1^2}-\displaystyle\frac{2\beta_2A_3}{\beta_1}
-\left(\displaystyle\frac{2\beta_3}{\beta_1}+\displaystyle\frac{\beta_2^2}{\beta_1^2}\right)A_4+\\
\qquad\qquad -2\left(\displaystyle\frac{\beta_4}{\beta_1}+\displaystyle\frac{\beta_2\beta_3}{\beta_1^2}\right)A_5, \\\\
A_1:=\displaystyle\frac{(b-c)^2-(a-c)^2}{2\beta_1}-\displaystyle\frac{\beta_2A_2}{\beta_1}
-\displaystyle\frac{\beta_3A_3}{\beta_1}-\displaystyle\frac{\beta_4A_4}{\beta_1}-\displaystyle\frac{\beta_5A_5}{\beta_1},
\end{array}\right.
\end{equation}  
then the 
$3$-tuple  $\big(\Lambda,\, \Omega_{\beta_1, \dots, \beta_5}, \Gamma_{\beta_1, \dots, \beta_5}\big)$ 
is the  $5$th-order automatic integration induced by $\mathcal{R}^{(6)}$, where the maps  $\Omega_{\beta_1, \dots, \beta_5}: \mathcal{R}\to \mathcal{R}^{(6)}$ and
$\Gamma_{\beta_1, \dots, \beta_5}: \mathcal{R}^{(6)}\to \mathcal{R}$ are defined by
\begin{equation}\label{eq6}
 \Omega_{\beta_1, \dots, \beta_5}(x):=x+\displaystyle\sum_{i=1}^5\beta_i\varepsilon^i
 \quad\mbox{for $x\in \mathcal{R}$}
\end{equation}
and 
\begin{equation}\label{eq7}
 \Gamma_{\beta_1, \dots, \beta_5}\left(\displaystyle\sum_{i=0}y_i\varepsilon^i\right):=(b-a)y_0+\sum_{i=1}^5 A_i y_i
 \quad\mbox{for $y_0$, $y_1$, $\dots$, $y_5\in \mathcal{R}$}.
\end{equation}
\end{proposition}

\medskip
\noindent
{\bf Proof}  First, by (\ref{eq5}) and  (\ref{eq6}),  the map $(\Lambda |{\bf D}^5(x)): {\bf D}^5(x)\to{\bf F}(\mathcal{R}^{(6)}_x, \,\mathcal{R}^{(6)})$  is clearly a real linear transformation. Hence,  the property (i)  in Definition \ref{def1.1}  holds.

\bigskip
Next,   for  $f$, $g\in {\bf D}^5(x)$,  we have
\begin{equation}\label{eq11}
\left\{\begin{array}{c}
(fg)'=f'g+fg',\\\\
(fg)''=f''g+2f'g'+fg'',\\\\
(fg)'''=f'''g+3f''g'+3f'g''+fg''',\\\\
(fg)^{(4)}=f^{(4)}g+4f^{(3)}g'+6f''g''+4f'g^{(3)}+fg^{(4)},\\\\
(fg)^{(5)}=f^{(5)}g+5f^{(4)}g'+10f^{(3)}g''+10f''g^{(3)}+5f'g^{(4)}+fg^{(5)}.
\end{array}\right.
\end{equation}

Let $x+\displaystyle\sum_{i=1}^5 a_i\varepsilon^i\in \mathcal{R}^{(6)}$, where 
$x$, $a_i\in \mathcal{R}$ for $1\le i\le 5$. By (\ref{eq5}),  we have
$$
\Big(\Lambda (f)\cdot \Lambda (g)\Big)\left(x+\displaystyle\sum_{i=1}^5 a_i\varepsilon^i \right)
=\Lambda (f)\left(x+\displaystyle\sum_{i=1}^5 a_i\varepsilon^i \right)\cdot  \Lambda (g)\left(x+\displaystyle\sum_{i=1}^5 a_i\varepsilon^i \right)
$$
\begin{eqnarray*}
&=&\left\{f+a_1f'\varepsilon+
+\Big(a_2f'+\displaystyle\frac{1}{2} a_1^2f^{(2)}\Big)\varepsilon^2+\Big(a_3f'+a_1a_2f^{(2)}+\displaystyle\frac{1}{6} a_1^3f^{(3)}\Big)\varepsilon^3+\right.\\
&&+\Big[a_4f'+\Big(a_1a_3+\displaystyle\frac{1}{2}a_2^2 \Big)f^{(2)}+
\displaystyle\frac{1}{2}a_1^2a_2 f^{(3)}+\displaystyle\frac{1}{24}a_1^4 f^{(4)}\Big]\varepsilon^4+
\end{eqnarray*}
$$
+\left.[a_5f'+(a_1a_4+a_2 a_3)f^{(2)}+
\displaystyle\frac{1}{2}(a_1^2a_3+a_1a_2^2) f^{(3)}+\displaystyle\frac{1}{6}a_1^3a_2 f^{(4)}+
\displaystyle\frac{1}{120}a_1^5 f^{(5)}]\varepsilon^5\right\}
$$
\begin{eqnarray*}
&\cdot&\left\{g+a_1g'\varepsilon+
+\Big(a_2g'+\displaystyle\frac{1}{2} a_1^2g^{(2)}\Big)\varepsilon^2+\Big(a_3g'+a_1a_2g^{(2)}+\displaystyle\frac{1}{6} a_1^3g^{(3)}\Big)\varepsilon^3+\right.\\
&&+\Big[a_4g'+\Big(a_1a_3+\displaystyle\frac{1}{2}a_2^2 \Big)g^{(2)}+
\displaystyle\frac{1}{2}a_1^2a_2 g^{(3)}+\displaystyle\frac{1}{24}a_1^4 g^{(4)}\Big]\varepsilon^4+
\end{eqnarray*}
$$
+\left.[a_5g'+(a_1a_4+a_2 a_3)g^{(2)}+
\displaystyle\frac{1}{2}(a_1^2a_3+a_1a_2^2) g^{(3)}+\displaystyle\frac{1}{6}a_1^3a_2 g^{(4)}+
\displaystyle\frac{1}{120}a_1^5 g^{(5)}]\varepsilon^5\right\}
$$
\begin{eqnarray*}
&=&fg+(a_1f'\cdot g+f\cdot a_1g')\varepsilon+\nonumber\\
&&\quad +\Big[f\cdot \Big(\underbrace{a_2g'}_{1}+\underbrace{\displaystyle\frac12 a_1^2g''}_{2}\Big)+ \underbrace{a_1f'\cdot a_1g'}_{2}
+\Big(\underbrace{a_2f'}_{1}+\underbrace{\displaystyle\frac12 a_1^2f''}_{2}\Big)\cdot g\Big]\varepsilon^2+\nonumber\\
&&\quad +\Big[f\cdot \Big(\underbrace{a_3g'}_{3}+\underbrace{a_1a_2g''}_{4}+
\underbrace{\displaystyle\frac16 a_1^3g'''}_{5}\Big)
+a_1f'\cdot \Big(\underbrace{a_2g'}_{4}+\underbrace{\displaystyle\frac12 a_1^2g''}_{5}\Big)+
\end{eqnarray*}
$$
+ \Big(\underbrace{a_2f'(x)}_{4}+\underbrace{\displaystyle\frac12 a_1^2f''}_{5}\Big)\cdot a_1g'+
\Big(\underbrace{a_3f'}_{3}+\underbrace{a_1a_2f''}_{4}+\underbrace{\displaystyle\frac16 a_1^3f'''}_{5}\Big)\cdot g\Big]\varepsilon^3+
$$
$$
+\left\{ f\cdot\Big[a_4g'+\Big(a_1a_3+\displaystyle\frac{1}{2}a_2^2 \Big)g^{(2)}+
\displaystyle\frac{1}{2}a_1^2a_2 g^{(3)}+\displaystyle\frac{1}{24}a_1^4 g^{(4)}\Big]+\right.
$$
$$
+a_1f'\cdot \Big(a_3g'+a_1a_2g^{(2)}+\displaystyle\frac{1}{6} a_1^3g^{(3)}\Big)+
\Big(a_2f'+\displaystyle\frac{1}{2} a_1^2f^{(2)}\Big)\cdot \Big(a_2g'+\displaystyle\frac{1}{2} a_1^2g^{(2)}\Big)+
$$
$$
+\Big(a_3f'+a_1a_2f^{(2)}+\displaystyle\frac{1}{6} a_1^3f^{(3)}\Big)\cdot a_1g'+
$$
$$
\left.\Big[a_4f'+\Big(a_1a_3+\displaystyle\frac{1}{2}a_2^2 \Big)f^{(2)}+
\displaystyle\frac{1}{2}a_1^2a_2 f^{(3)}+\displaystyle\frac{1}{24}a_1^4 f^{(4)}\Big]\cdot g\right\}
\varepsilon^4+
$$
$$
+\left\{f \Big[a_5g'+(a_1a_4+a_2 a_3)g^{(2)}+
\displaystyle\frac{1}{2}(a_1^2a_3+a_1a_2^2) g^{(3)}+\displaystyle\frac{1}{6}a_1^3a_2 g^{(4)}+
\displaystyle\frac{1}{120}a_1^5 g^{(5)}\Big]\right.+
$$
$$
+a_1f'\cdot \Big[a_4g'+\Big(a_1a_3+\displaystyle\frac{1}{2}a_2^2 \Big)g^{(2)}+
\displaystyle\frac{1}{2}a_1^2a_2 g^{(3)}+\displaystyle\frac{1}{24}a_1^4 g^{(4)}\Big]+
$$
$$
+\Big(a_2f'+\displaystyle\frac{1}{2} a_1^2f^{(2)}\Big)\cdot 
\Big(a_3g'+a_1a_2g^{(2)}+\displaystyle\frac{1}{6} a_1^3g^{(3)}\Big)+
$$
$$
+\Big(a_3f'+a_1a_2f^{(2)}+\displaystyle\frac{1}{6} a_1^3f^{(3)}\Big)\cdot 
\Big(a_2g'+\displaystyle\frac{1}{2} a_1^2g^{(2)}\Big)+
$$
$$
+\Big[a_4f'+\Big(a_1a_3+\displaystyle\frac{1}{2}a_2^2 \Big)f^{(2)}+
\displaystyle\frac{1}{2}a_1^2a_2 f^{(3)}+\displaystyle\frac{1}{24}a_1^4 f^{(4)}\Big]\cdot a_1g'+
$$
\begin{eqnarray}\label{eq12}
&&\quad +\Big[a_5f'+(a_1a_4+a_2 a_3)f^{(2)}+
\displaystyle\frac{1}{2}(a_1^2a_3+a_1a_2^2) f^{(3)}+\nonumber\\
&&\quad\qquad +\left.\displaystyle\frac{1}{6}a_1^3a_2 f^{(4)}+
\displaystyle\frac{1}{120}a_1^5 f^{(5)}\Big] g\right\}\varepsilon^5.
\end{eqnarray}

\medskip
By (\ref{eq11}), we have
\begin{eqnarray}\label{eq13}
\left\{\begin{array}{l}
\mbox{the coefficient of $\varepsilon$ in (\ref{eq12})} =a_1(fg)',\\
\mbox{the coefficient of $\varepsilon^2$  in (\ref{eq12})}=a_2(fg)'+\displaystyle\frac12 a_1^2(fg)'',\\
\mbox{the coefficient of $\varepsilon^3$  in (\ref{eq12})}=a_3(fg)'+a_1a_2(fg)''+\displaystyle\frac16 a_1^3(fg)''',
\end{array}\right.
\end{eqnarray}
\begin{eqnarray}\label{eq14}
&&\mbox{the coefficient of $\varepsilon^4$  in (\ref{eq12})}=\nonumber\\
&=&a_4(fg)'+\Big(a_1a_3+\displaystyle\frac{1}{2}a_2^2 \Big)(fg)'' 
+\displaystyle\frac{1}{2}a_1^2a_2(fg)^{(3)}+\displaystyle\frac{1}{24}a_1^4(fg)^{(4)},
\end{eqnarray}
\begin{eqnarray}\label{eq16}
&&\mbox{the coefficient of $\varepsilon^5$  in (\ref{eq12})}
=a_5(fg)'+(a_1a_4+a_2 a_3)(fg)'' +\nonumber\\
&&\quad +\displaystyle\frac{1}{2}(a_1^2a_3+a_1a_2^2)(fg)^{(3)}
+\displaystyle\frac{1}{6}a_1^3a_2(fg)^{(4)}
 +\displaystyle\frac{1}{120}a_1^5(fg)^{(5)}.
\end{eqnarray}

\medskip
For example, let us check  (\ref{eq16}). Using  (\ref{eq12}), we get
\begin{eqnarray}\label{eq15}
&&\mbox{the coefficient of $\varepsilon^5$  in (\ref{eq12})}=\underbrace{a_5fg'}_1+\underbrace{(a_1a_4+a_2 a_3)fg^{(2)}}_2+\nonumber\\
&&\quad +
\underbrace{\displaystyle\frac{1}{2}(a_1^2a_3+a_1a_2^2) fg^{(3)}}_3+\underbrace{\displaystyle\frac{1}{6}a_1^3a_2f g^{(4)}}_4+
\underbrace{\displaystyle\frac{1}{120}a_1^5f g^{(5)}}_5+
\underbrace{a_1a_4f'g'}_2+\nonumber\\
&&\quad +\underbrace{\Big(a_1^2a_3+\displaystyle\frac{1}{2}a_1a_2^2 \Big)f'g^{(2)}}_3+
\underbrace{\displaystyle\frac{1}{2}a_1^3a_2 f'g^{(3)}}_4+\underbrace{\displaystyle\frac{1}{24}a_1^5f' g^{(4)}}_5
+\underbrace{a_2a_3f'g'}_2+\nonumber\\
&&\quad+\underbrace{a_1a_2^2f'g^{(2)}}_3+\underbrace{\displaystyle\frac{1}{6} a_1^3a_2f'g^{(3)}}_4+
\underbrace{\displaystyle\frac{1}{2} a_1^2a_3f^{(2)}g'}_3+\underbrace{\displaystyle\frac{1}{2}a_1^3a_2f^{(2)}g^{(2)}}_4+\nonumber\\
&&\quad +\underbrace{\displaystyle\frac{1}{12} a_1^5f^{(2)}g^{(3)}}_5+
\underbrace{a_2a_3f'g'}_2+\underbrace{\displaystyle\frac{1}{2} a_1^2a_3f'g^{(2)}}_3+
\underbrace{a_1a_2^2f^{(2)}g'}_3+\nonumber\\
&&\quad +\underbrace{\displaystyle\frac{1}{2} a_1^3a_2f^{(2)}g^{(2)}}_4+\underbrace{\displaystyle\frac{1}{6}a_1^3a_2f^{(3)}g'}_4+\underbrace{\displaystyle\frac{1}{12} a_1^5f^{(3)}g^{(2)}}_5+
\underbrace{a_1a_4f'g'}_2+\nonumber\\
&&\quad +\underbrace{\Big(a_1^2a_3+\displaystyle\frac{1}{2}a_1a_2^2 \Big)f^{(2)}g'}_3+
\underbrace{\displaystyle\frac{1}{2}a_1^3a_2 f^{(3)}g'}_4+\underbrace{\displaystyle\frac{1}{24}a_1^5 f^{(4)}g'}_5+
\underbrace{a_5f'g}_1+\nonumber\\
&& +\underbrace{(a_1a_4+a_2 a_3)f^{(2)}g}_2+
\underbrace{\displaystyle\frac{1}{2}(a_1^2a_3+a_1a_2^2) f^{(3)}g}_3+\underbrace{\displaystyle\frac{1}{6}a_1^3a_2 f^{(4)}g}_4+
\underbrace{\displaystyle\frac{1}{120}a_1^5 f^{(5)}g}_5\nonumber\\
&=&a_5(f'g+fg')+(a_1a_4+a_2 a_3)(f''g+2f'g'+fg'')+\nonumber\\
&&\quad +\displaystyle\frac{1}{2}(a_1^2a_3+a_1a_2^2)(f'''g+3f''g'+3f'g''+fg''')+\nonumber\\
&&\quad\quad +\displaystyle\frac{1}{6}a_1^3a_2\Big(f^{(4)}g+4f^{(3)}g'+6f''g''+4f'g^{(3)}+fg^{(4)}\Big)+\nonumber\\
&& +\displaystyle\frac{1}{120}a_1^5\Big(f^{(5)}g+5f^{(4)}g'+10f^{(3)}g''+10f''g^{(3)}+5f'g^{(4)}+fg^{(5)}\Big),
\end{eqnarray}
which prove sthat  (\ref{eq16}) holds.

\medskip
By (\ref{eq5}),  (\ref{eq12}),  (\ref{eq13}),  (\ref{eq14}) and  (\ref{eq16}),  we get
\begin{equation}\label{eq17}
\Lambda (f)\cdot \Lambda (g)=\Lambda (f\cdot g) \quad\mbox{for $f$, $g\in   {\bf D}^5(x)$,}
\end{equation}
which proves  that  $\Lambda$  preserves the product  in the algebra  ${\bf D}^5(x)$. In particular, the property (ii)  in Definition \ref{def1.1}  holds.

\bigskip
Thirdly, let  $x\in \mathcal{R}$,  $g\in {\bf D}^5x)$ and $f\in {\bf D}^5\big(g(x)\big)$. It follows from  (\ref{eq5}) and (\ref{eq6}) that
$$
\Big(\Lambda (f)\circ \Lambda (g)\Big)\left(x+\displaystyle\sum_{i=1}^5 a_i\varepsilon^i \right)
=\Lambda (f)\left( \Lambda (g)\left(x+\displaystyle\sum_{i=1}^5 a_i\varepsilon^i \right)\right)
$$
$$
=\Lambda (f)\left(g+a_1g'\varepsilon+ \Big(a_2g'+\displaystyle\frac{1}{2} a_1^2g^{(2)}\Big)\varepsilon^2+\Big(a_3g'+a_1a_2g^{(2)}+\displaystyle\frac{1}{6} a_1^3g^{(3)}\Big)\varepsilon^3\right.+
$$
$$
+\Big[a_4g'+\Big(a_1a_3+\displaystyle\frac{1}{2}a_2^2 \Big)g^{(2)}+
\displaystyle\frac{1}{2}a_1^2a_2 g^{(3)}+\displaystyle\frac{1}{24}a_1^4 g^{(4)}\Big]\varepsilon^4+
$$
$$
+\left.\Big[a_5g'+(a_1a_4+a_2 a_3)g^{(2)}+
\displaystyle\frac{1}{2}(a_1^2a_3+a_1a_2^2) g^{(3)}+\displaystyle\frac{1}{6}a_1^3a_2 g^{(4)}+
\displaystyle\frac{1}{5!}a_1^5 g^{(5)}\Big]\varepsilon^5\right)
$$
\begin{eqnarray*}
&=&f(g)+a_1g'f'(g)'\varepsilon+\Big\{\Big(a_2g'+\displaystyle\frac{1}{2} a_1^2g^{(2)}\Big)f'(g)+
\displaystyle\frac{1}{2}(a_1g')^2f''(g)\Big\}\varepsilon^2+\\
&& +\Big\{\Big(a_3g'+a_1a_2g^{(2)}+\displaystyle\frac{1}{6} a_1^3g^{(3)}\Big)\cdot f'(g)+
a_1g'\cdot \Big(a_2g'+\displaystyle\frac{1}{2} a_1^2g^{(2)}\Big)\cdot f''(g)+\\
&&\quad +\displaystyle\frac{1}{6}(a_1g')^3f^{(3)}(g)\Big\}\varepsilon^3+
\Big\{a_4g'+\Big(a_1a_3+\displaystyle\frac{1}{2}a_2^2 \Big)g^{(2)}+
\displaystyle\frac{1}{2}a_1^2a_2 g^{(3)}+\\
&& +\displaystyle\frac{1}{24}a_1^4 g^{(4)}\Big]f'(g)+\Big[a_1g'\cdot \Big(a_3g'+a_1a_2g^{(2)}+\displaystyle\frac{1}{6} a_1^3g^{(3)}\Big)+\\
&&\quad +\displaystyle\frac{1}{2} \Big(a_2g'+\displaystyle\frac{1}{2} a_1^2g^{(2)}\Big)^2\Big]f''(g)+
\displaystyle\frac{1}{2}(a_1g')^2\cdot \Big(a_2g'+\displaystyle\frac{1}{2} a_1^2g^{(2)}\Big)f^{(3)}(g)+\\
&& +\displaystyle\frac{1}{24}(a_1g')^4f^{(4)}(g)\Big\}\varepsilon^4+
\Big\{\Big[a_5g'+(a_1a_4+a_2 a_3)g^{(2)}+\\
&&\quad +\displaystyle\frac{1}{2}(a_1^2a_3+a_1a_2^2) g^{(3)}+
\displaystyle\frac{1}{6}a_1^3a_2 g^{(4)}+
\displaystyle\frac{1}{120}a_1^5 g^{(5)}\Big]f'(g)+\\
&& +\Big[a_1g'\cdot\Big(a_4g'+\Big(a_1a_3+\displaystyle\frac{1}{2}a_2^2 \Big)g^{(2)}+
\displaystyle\frac{1}{2}a_1^2a_2 g^{(3)}+\displaystyle\frac{1}{24}a_1^4 g^{(4)}\Big)+\\
&&\quad +\Big(a_2g'+\displaystyle\frac{1}{2} a_1^2g^{(2)}\Big)\cdot \Big(a_3g'+a_1a_2g^{(2)}+\displaystyle\frac{1}{6} a_1^3g^{(3)}\Big)\Big]f''(g)+
\end{eqnarray*}
$$
+\displaystyle\frac{1}{2}\Big[(a_1g')^2\cdot \Big(a_3g'+a_1a_2g^{(2)}+\displaystyle\frac{1}{6} a_1^3g^{(3)}\Big)+a_1g'\cdot  \Big(a_2g'+\displaystyle\frac{1}{2} a_1^2g^{(2)}\Big)^2 \Big]f^{(3)}(g)+
$$
\begin{equation}\label{eq18}
+\displaystyle\frac{1}{6}(a_1g')^3\cdot  \Big(a_2g'+\displaystyle\frac{1}{2} a_1^2g^{(2)}\Big)f^{(4)}(g)+
\displaystyle\frac{1}{120}(a_1g')^5f^{(5)}(g)\Big\}\varepsilon^5.
\end{equation}

\medskip
By the chain rule, we have
\begin{eqnarray}
(f\circ g)'&=&g'f'(g),\nonumber\\
(f\circ g)''&=&g''f'(g)+(g')^2f''(g),,\nonumber\\
(f\circ g)^{(3)}&=&g^{(3)}f'(g)+3g'g''f''(g)+(g')^3f^{(3)}(g),,\nonumber\\
(f\circ g)^{(4)}&=&g^{(4)}f'(g)+4g'g^{(3)}f''(g)+3(g'')^2f''(g)+,\nonumber\\
&&\qquad\qquad\qquad +6(g')^2g''f^{(3)}(g)+(g')^4f^{(4)}(g),\label{eq19}\\
(f\circ g)^{(5)}&=&g^{(5)}f'(g)+5g'g^{(4)}f''(g)+10g''g^{(3)}f''(g)+
10(g')^2g^{(3)}f^{(3)}(g)+\nonumber\\
&&\qquad+15g'(g'')^2f^{(3)}(g) +10(g')^3g''f^{(4)}(g)+(g')^5f^{(5)}(g).\label{eq20}
\end{eqnarray}

\medskip
It follows from  the facts above and (\ref{eq18}) that
\begin{eqnarray}\label{eq21}
\left\{\begin{array}{l}
\mbox{the coefficient of $\varepsilon$ in (\ref{eq18})} =a_1(f\circ g)',\\
\mbox{the coefficient of $\varepsilon^2$  in (\ref{eq18})}=a_2(f\circ g)'+\displaystyle\frac12 a_1^2(f\circ g)'',\\
\mbox{the coefficient of $\varepsilon^3$  in (\ref{eq18})}=a_3(f\circ g)'+a_1a_2(f\circ g)''+\displaystyle\frac16 a_1^3(f\circ g)''',
\end{array}\right.
\end{eqnarray}
\begin{eqnarray}\label{eq22}
&&\mbox{the coefficient of $\varepsilon^4$  in (\ref{eq18})}
=a_4(f\circ g)'+\nonumber\\
&&\quad +\Big(a_1a_3+\displaystyle\frac{1}{2}a_2^2 \Big)(f\circ g)''+\displaystyle\frac{1}{2}a_1^2a_2(f\circ g)^{(3)}+
\displaystyle\frac{1}{24}a_1^4(f\circ g)^{(4)},
\end{eqnarray}

\begin{eqnarray}\label{eq23}
&&\mbox{the coefficient of $\varepsilon^5$  in (\ref{eq18})}=a_5(f\circ g)'+(a_1a_4+a_2 a_3)(f\circ g)''+
\nonumber\\
&&\quad +\displaystyle\frac{1}{2}(a_1^2a_3+a_1a_2^2)(f\circ g)^{(3)}+
\displaystyle\frac{1}{6}a_1^3a_2(f\circ g)^{(4)}+\displaystyle\frac{1}{120}a_1^5(f\circ g)^{(5)}.
\end{eqnarray}

Using  (\ref{eq21}),  (\ref{eq22}) and (\ref{eq23}), the equation  (\ref{eq18}) becomes
\begin{eqnarray*}
&&\Big(\Lambda (f)\circ \Lambda (g)\Big)\left(x+\displaystyle\sum_{i=1}^5 a_i\varepsilon^i \right)\\
&=&f(g)+a_1(f\circ g)'\varepsilon+\Big\{a_2(f\circ g)'+\displaystyle\frac12 a_1^2(f\circ g)''\Big\}\varepsilon^2+\\
&&\quad +\Big\{a_3(f\circ g)'+a_1a_2(f\circ g)''+\displaystyle\frac16 a_1^3(f\circ g)'''\Big\}\varepsilon^3+\\
&&+\Big\{a_4(f\circ g)'+\Big(a_1a_3+\displaystyle\frac{1}{2}a_2^2 \Big)(f\circ g)''+\displaystyle\frac{1}{2}a_1^2a_2(f\circ g)^{(3)}+\\
&&\quad +\displaystyle\frac{1}{24}a_1^4(f\circ g)^{(4)}\Big\}\varepsilon^4+\Big\{a_5(f\circ g)'+(a_1a_4+a_2 a_3)(f\circ g)''+\\
&& +
\displaystyle\frac{1}{2}(a_1^2a_3+a_1a_2^2)(f\circ g)^{(3)}+
\displaystyle\frac{1}{6}a_1^3a_2(f\circ g)^{(4)}+\displaystyle\frac{1}{120}a_1^5(f\circ g)^{(5)}\Big\}\varepsilon^5\\
&=&\Big(\Lambda (f\circ g)\Big)\left(x+\displaystyle\sum_{i=1}^5 a_i\varepsilon^i \right),
\end{eqnarray*}
which implies that  the property (iii)  in Definition \ref{def1.1} holds.

\bigskip
Finally, let 
$$
T_{5; c}(x):=f(c)+f'(c)(x-c)+\displaystyle\frac{f^{(2)}(c)}{2!}(x-c)^{2}+\cdots +
\displaystyle\frac{f^{(5)}(c)}{5!}(x-c)^{5}
$$
be the  $5$th-order Taylor polynomial for $f(x)$ with its  center at $c$. Then we have
\begin{eqnarray*}
\displaystyle\int T_{5; c}(x)dx&=&\displaystyle\int \left(\displaystyle\sum_{i=0}^5\displaystyle\frac{f^{(i)}(c)}{i!}(x-c)^{i}\right)dx=\displaystyle\sum_{i=0}^5\displaystyle\int\displaystyle\frac{f^{(i)}(c)}{i!}(x-c)^{i}dx\\
&=&K+\displaystyle\sum_{i=0}^5\displaystyle\frac{f^{(i)}(c)}{(i+1)!}(x-c)^{i+1}\quad\mbox{ for some constany $K$},
\end{eqnarray*}
which implies that $\displaystyle\int_a^b T_{5; c}(x)dx=
\left.\left(K+\displaystyle\sum_{i=0}^5\displaystyle\frac{f^{(i)}(c)}{(i+1)!}(x-c)^{i+1}\right)\right|_a^b$ or
\begin{equation}\label{eq24}
\displaystyle\int_a^b T_{5; c}(x)dx
=\displaystyle\sum_{i=0}^5\displaystyle\frac{(b-c)^{i+1}-(a-c)^{i+1}}{(i+1)!}f^{(i)}(c).
\end{equation}

\medskip
By  (\ref{eq6}) and (\ref{eq7}), we get
\begin{eqnarray}\label{eq25}
&&\Big(\Lambda(f)\circ\Omega_{\beta_1, \dots, \beta_5}\Big)(c)=\Lambda(f)\left(c+\displaystyle\sum_{i=1}^5\beta_i\varepsilon^i\right)\nonumber\\
&=&\underbrace{f(c)}_{y_0}+\underbrace{\beta_1f'(c)}_{y_1}\varepsilon+
+\underbrace{\Big\{\beta_2f'(c)+\displaystyle\frac{1}{2!} \beta_1^2f^{(2)}(c)\Big\}}_{y_2}\varepsilon^2+\underbrace{\Big\{\beta_3f'(c)+\beta_1\beta_2f^{(2)}(c)}_{y_3}+\nonumber\\
&&\quad +\underbrace{\displaystyle\frac{1}{3!} \beta_1^3f^{(3)}(c)\Big\}}_{y_3}\varepsilon^3
+\underbrace{\Big\{\beta_4f'(c)+\Big(\beta_1\beta_3+\displaystyle\frac{1}{2}\beta_2^2 \Big)f^{(2)}(c)+
\displaystyle\frac{1}{2}\beta_1^2\beta_2 f^{(3)}(c)}_{y_4}+\nonumber\\
&&\quad\quad +\underbrace{\displaystyle\frac{1}{4!}\beta_1^4 f^{(4)}(c)\Big\}}_{y_4}\varepsilon^4
+\underbrace{\Big\{\beta_5f'(c)+(\beta_1\beta_4+\beta_2 \beta_3)f^{(2)}(c)}_{y_5}+\nonumber\\
&&\quad +
\underbrace{\displaystyle\frac{1}{2}(\beta_1^2\beta_3+\beta_1\beta_2^2) f^{(3)}(c)+\displaystyle\frac{1}{6}\beta_1^3\beta_2 f^{(4)}(c)+
\displaystyle\frac{1}{5!}\beta_1^5 f^{(5)}(c)\Big\}}_{y_5}\varepsilon^5.
\end{eqnarray}

It follows from  (\ref{eq8}) and (\ref{eq25}) that
\begin{eqnarray}\label{eq26}
&&\Big(\Gamma_{\beta_1, \dots, \beta_5}\circ\Lambda(f)\circ\Omega_{\beta_1, \dots, \beta_5}\Big)(c)=
(b-a)y_0+ \displaystyle\sum_{i=1}^5 A_i y_i\nonumber\\
&=&(b-a)\cdot f(c)+\left[\displaystyle\frac{(b-c)^2-(a-c)^2}{2\beta_1}-\displaystyle\frac{\beta_2A_2}{\beta_1}
-\displaystyle\frac{\beta_3A_3}{\beta_1}-\displaystyle\frac{\beta_4A_4}{\beta_1}+\right.\nonumber\\
&&\quad \left.-\displaystyle\frac{\beta_5A_5}{\beta_1}\right]\cdot \beta_1\underbrace{f'(c)}_1+\left[\displaystyle\frac{(b-c)^3-(a-c)^3}{3\beta_1^2}-\displaystyle\frac{2\beta_2A_3}{\beta_1}\right.+\nonumber\\
&& -\left.\left(\displaystyle\frac{2\beta_3}{\beta_1}+\displaystyle\frac{\beta_2^2}{\beta_1^2}\right)A_4-2\left(\displaystyle\frac{\beta_4}{\beta_1}+\displaystyle\frac{\beta_2\beta_3}{\beta_1^2}\right)A_5\right]\cdot \Big\{\beta_2\underbrace{f'(c)}_1+\displaystyle\frac{1}{2!} \beta_1^2\underbrace{f^{(2)}(c)}_2\Big\}+\nonumber\\
&& \quad +\left[\displaystyle\frac{(b-c)^4-(a-c)^4}{4\beta_1^3}-\displaystyle\frac{3\beta_2A_4}{\beta_1}
-3\left(\displaystyle\frac{\beta_3}{\beta_1}+\displaystyle\frac{\beta_2^2}{\beta_1^2}\right)A_5\right]\cdot
\Big\{\beta_3\underbrace{f'(c)}_1+\nonumber\\
&&+\beta_1\beta_2\underbrace{f^{(2)}(c)}_2+\displaystyle\frac{1}{3!} \beta_1^3\underbrace{f^{(3)}(c)}_3\Big\}+
\left[\displaystyle\frac{(b-c)^5-(a-c)^5}{5\beta_1^4}-\displaystyle\frac{4\beta_2A_5}{\beta_1}\right]\cdot
\Big\{\beta_4\underbrace{f'(c)}_1+\nonumber\\
&&\quad +\Big(\beta_1\beta_3+\displaystyle\frac{1}{2}\beta_2^2 \Big)\underbrace{f^{(2)}(c)}_2+
\displaystyle\frac{1}{2}\beta_1^2\beta_2 \underbrace{f^{(3)}(c)}_3 +\displaystyle\frac{1}{4!}\beta_1^4 \underbrace{f^{(4)}(c)}_4\Big\}+\nonumber\\
&&+
\left[\displaystyle\frac{(b-c)^6-(a-c)^6}{6\beta_1^5}\right]\cdot \Big\{\beta_5\underbrace{f'(c)}_1+(\beta_1\beta_4+\beta_2 \beta_3)\underbrace{f^{(2)}(c)}_2+\nonumber\\
&&\qquad +\displaystyle\frac{1}{2}(\beta_1^2\beta_3+\beta_1\beta_2^2) \underbrace{f^{(3)}(c)}_3+\displaystyle\frac{1}{6}\beta_1^3\beta_2\underbrace{ f^{(4)}(c)}_4+
\displaystyle\frac{1}{5!}\beta_1^5\underbrace{ f^{(5)}(c)}_5\Big\}.
\end{eqnarray}

Using  (\ref{eq26}), we have
\begin{eqnarray*}
&&\mbox{the coefficient of $f'(c)$}=\left[\displaystyle\frac{(b-c)^2-(a-c)^2}{2\beta_1}-\displaystyle\frac{\beta_2A_2}{\beta_1}
-\displaystyle\frac{\beta_3A_3}{\beta_1}-\displaystyle\frac{\beta_4A_4}{\beta_1}+\right.\nonumber\\
&&\quad \left.-\displaystyle\frac{\beta_5A_5}{\beta_1}\right]\cdot \beta_1+
\left[\displaystyle\frac{(b-c)^3-(a-c)^3}{3\beta_1^2}-\displaystyle\frac{2\beta_2A_3}{\beta_1} -\left(\displaystyle\frac{2\beta_3}{\beta_1}+\displaystyle\frac{\beta_2^2}{\beta_1^2}\right)A_4\right.+\\
&&\qquad \left.-2\left(\displaystyle\frac{\beta_4}{\beta_1}+\displaystyle\frac{\beta_2\beta_3}{\beta_1^2}\right)A_5\right]\cdot \beta_2+
\left[\displaystyle\frac{(b-c)^4-(a-c)^4}{4\beta_1^3}-\displaystyle\frac{3\beta_2A_4}{\beta_1}+\right.\\
&&\quad -\left.3\left(\displaystyle\frac{\beta_3}{\beta_1}+\displaystyle\frac{\beta_2^2}{\beta_1^2}\right)A_5\right]\cdot\beta_3+\left[\displaystyle\frac{(b-c)^5-(a-c)^5}{5\beta_1^4}-\displaystyle\frac{4\beta_2A_5}{\beta_1}\right]\cdot\beta_4+\\
&&\quad +\displaystyle\frac{(b-c)^6-(a-c)^6}{6\beta_1^5}\cdot \beta_5\\
&=&
\displaystyle\frac{(b-c)^2-(a-c)^2}{2}-\beta_2A_2
-\beta_3A_3+\\
&&\quad -\beta_4A_4 -\beta_5A_5+
\left[\displaystyle\frac{(b-c)^3-(a-c)^3}{3\beta_1^2}-\displaystyle\frac{2\beta_2A_3}{\beta_1} -\left(\displaystyle\frac{2\beta_3}{\beta_1}+\displaystyle\frac{\beta_2^2}{\beta_1^2}\right)A_4\right.+\\
&&\qquad \left.-2\left(\displaystyle\frac{\beta_4}{\beta_1}+\displaystyle\frac{\beta_2\beta_3}{\beta_1^2}\right)A_5\right]\cdot \beta_2+
\left[\displaystyle\frac{(b-c)^4-(a-c)^4}{4\beta_1^3}-\displaystyle\frac{3\beta_2A_4}{\beta_1}+\right.\\
&&\quad -\left.3\left(\displaystyle\frac{\beta_3}{\beta_1}+\displaystyle\frac{\beta_2^2}{\beta_1^2}\right)A_5\right]\cdot\beta_3+\left[\displaystyle\frac{(b-c)^5-(a-c)^5}{5\beta_1^4}-\displaystyle\frac{4\beta_2A_5}{\beta_1}\right]\cdot\beta_4+\\
&&\quad +\displaystyle\frac{(b-c)^6-(a-c)^6}{6\beta_1^5}\cdot \beta_5\\
&\stackrel{  (\ref{eq8})}{=}&
\displaystyle\frac{(b-c)^2-(a-c)^2}{2}-
\beta_2\left\{\underbrace{\displaystyle\frac{(b-c)^3-(a-c)^3}{3\beta_1^2}-
\displaystyle\frac{2\beta_2A_3}{\beta_1}}_{1}\right.+\\
&& -\left.\underbrace{\left(\displaystyle\frac{2\beta_3}{\beta_1}
+\displaystyle\frac{\beta_2^2}{\beta_1^2}\right)A_4 -2\left(\displaystyle\frac{\beta_4}{\beta_1}
+\displaystyle\frac{\beta_2\beta_3}{\beta_1^2}\right)A_5}_{1}\right\}-
\beta_3\left\{\underbrace{\displaystyle\frac{(b-c)^4-(a-c)^4}{4\beta_1^3}}_{2}+\right.\\
&&\left.-\underbrace{\displaystyle\frac{3\beta_2A_4}{\beta_1}
-3\left(\displaystyle\frac{\beta_3}{\beta_1}+\displaystyle\frac{\beta_2^2}{\beta_1^2}\right)A_5}_{2}\right\}
 -\beta_4\left\{\underbrace{\displaystyle\frac{(b-c)^5-(a-c)^5}{5\beta_1^4}-\displaystyle\frac{4\beta_2A_5}{\beta_1}}_{3}\right\}+\\
 &&\quad -\underbrace{\beta_5\cdot\displaystyle\frac{(b-c)^6-(a-c)^6}{6\beta_1^5}}_{4}+
 \left[\underbrace{\displaystyle\frac{(b-c)^3-(a-c)^3}{3\beta_1^2}-
 \displaystyle\frac{2\beta_2A_3}{\beta_1}}_{1}+\right.\\
 && \left.-\underbrace{\left(\displaystyle\frac{2\beta_3}{\beta_1}+\displaystyle\frac{\beta_2^2}{\beta_1^2}\right)A_4-2\left(\displaystyle\frac{\beta_4}{\beta_1}+\displaystyle\frac{\beta_2\beta_3}{\beta_1^2}\right)A_5}_{1}\right]\cdot \beta_2+
\left[\underbrace{\displaystyle\frac{(b-c)^4-(a-c)^4}{4\beta_1^3}}_{2}+\right.
\end{eqnarray*}
$$
\left.-\underbrace{\displaystyle\frac{3\beta_2A_4}{\beta_1}
 -3\left(\displaystyle\frac{\beta_3}{\beta_1}+\displaystyle\frac{\beta_2^2}{\beta_1^2}\right)A_5}_{2}\right]\cdot\beta_3+\left[\underbrace{\displaystyle\frac{(b-c)^5-(a-c)^5}{5\beta_1^4}-\displaystyle\frac{4\beta_2A_5}{\beta_1}}_{3}\right]\cdot\beta_4+
 $$
 \begin{equation}\label{eq27}
+\underbrace{\displaystyle\frac{(b-c)^6-(a-c)^6}{6\beta_1^5}\cdot \beta_5}_{4}
=\displaystyle\frac{(b-c)^2-(a-c)^2}{2},\qquad\qquad
\end{equation}
\begin{eqnarray*}
&&\mbox{the coefficient of $f''(c)$}\\
&=&
\left[\displaystyle\frac{(b-c)^3-(a-c)^3}{3\beta_1^2}-\displaystyle\frac{2\beta_2A_3}{\beta_1}+
 -\left(\displaystyle\frac{2\beta_3}{\beta_1}+\displaystyle\frac{\beta_2^2}{\beta_1^2}\right)A_4+\right.\\
 &&\quad \left.-2\left(\displaystyle\frac{\beta_4}{\beta_1}+\displaystyle\frac{\beta_2\beta_3}{\beta_1^2}\right)A_5\right]\cdot \displaystyle\frac{1}{2} \beta_1^2+\\
&&\qquad +\left[\displaystyle\frac{(b-c)^4-(a-c)^4}{4\beta_1^3}
 -\displaystyle\frac{3\beta_2A_4}{\beta_1}
-3\left(\displaystyle\frac{\beta_3}{\beta_1}+\displaystyle\frac{\beta_2^2}{\beta_1^2}\right)A_5\right]\cdot
\beta_1\beta_2+\\
&&\quad +
\left[\displaystyle\frac{(b-c)^5-(a-c)^5}{5\beta_1^4}-\displaystyle\frac{4\beta_2A_5}{\beta_1}\right]\cdot
\Big(\beta_1\beta_3+\displaystyle\frac{1}{2}\beta_2^2 \Big)+\\
&&\qquad +\displaystyle\frac{(b-c)^6-(a-c)^6}{6\beta_1^5}\cdot (\beta_1\beta_4+\beta_2 \beta_3)
\end{eqnarray*}
$$
=\displaystyle\frac{(b-c)^3-(a-c)^3}{6}-\beta_1\beta_2A_3
 -\left(\beta_1\beta_3+\displaystyle\frac{\beta_2^2}{2}\right)A_4
 -\left(\beta_1\beta_4+\beta_2\beta_3\right)A_5+
 $$
 \begin{eqnarray}\label{eq28}
&&\qquad +\underbrace{\left[\displaystyle\frac{(b-c)^4-(a-c)^4}{4\beta_1^3}
 -\displaystyle\frac{3\beta_2A_4}{\beta_1}
-3\left(\displaystyle\frac{\beta_3}{\beta_1}+\displaystyle\frac{\beta_2^2}{\beta_1^2}\right)A_5
\right]}_{A_3}\cdot
\beta_1\beta_2+\nonumber\\
&&\quad +
\underbrace{\left[\displaystyle\frac{(b-c)^5-(a-c)^5}{5\beta_1^4}-\displaystyle\frac{4\beta_2A_5}{\beta_1}\right]}_{A_4}\cdot
\Big(\beta_1\beta_3+\displaystyle\frac{1}{2}\beta_2^2 \Big)+\nonumber\\
&&\qquad +\underbrace{\displaystyle\frac{(b-c)^6-(a-c)^6}{6\beta_1^5}}_{A_5}\cdot (\beta_1\beta_4+\beta_2 \beta_3)=\displaystyle\frac{(b-c)^3-(a-c)^3}{3!},
\end{eqnarray}
\begin{eqnarray*}
&&\mbox{the coefficient of $f^{(3)}(c)$}\\
&=&\left[\displaystyle\frac{(b-c)^4-(a-c)^4}{4\beta_1^3}-\displaystyle\frac{3\beta_2A_4}{\beta_1}
-3\left(\displaystyle\frac{\beta_3}{\beta_1}+\displaystyle\frac{\beta_2^2}{\beta_1^2}\right)A_5\right]\cdot
\displaystyle\frac{1}{3!} \beta_1^3+
\end{eqnarray*}
$$
+\left[\displaystyle\frac{(b-c)^5-(a-c)^5}{5\beta_1^4}-\displaystyle\frac{4\beta_2A_5}{\beta_1}\right]
\displaystyle\frac{1}{2}\beta_1^2\beta_2 +\displaystyle\frac{(b-c)^6-(a-c)^6}{6\beta_1^5}\displaystyle\frac{1}{2}(\beta_1^2\beta_3+\beta_1\beta_2^2)
$$
\begin{eqnarray}\label{eq29}
&=&\displaystyle\frac{(b-c)^4-(a-c)^4}{4!}-\displaystyle\frac12 \beta_1\beta_2A_4
-\frac{1}{2}(\beta_1^2\beta_3+\beta_1\beta_2^2)A_5+\nonumber\\
&&\quad +\underbrace{\left[\displaystyle\frac{(b-c)^5-(a-c)^5}{5\beta_1^4}-\displaystyle\frac{4\beta_2A_5}{\beta_1}\right]}_{A_4}
\displaystyle\frac{1}{2}\beta_1^2\beta_2 +\nonumber\\
&& +
\underbrace{\displaystyle\frac{(b-c)^6-(a-c)^6}{6\beta_1^5}}_{A_5}\displaystyle\frac{1}{2}(\beta_1^2\beta_3+\beta_1\beta_2^2)=\displaystyle\frac{(b-c)^4-(a-c)^4}{4!},
\end{eqnarray}
\begin{eqnarray}\label{eq30}
&&\mbox{the coefficient of $f^{(4)}(c)$}\nonumber\\
&=&
\left[\displaystyle\frac{(b-c)^5-(a-c)^5}{5\beta_1^4}-\displaystyle\frac{4\beta_2A_5}{\beta_1}\right]\cdot
\displaystyle\frac{1}{4!}\beta_1^4 +
\displaystyle\frac{(b-c)^6-(a-c)^6}{6\beta_1^5}\cdot\displaystyle\frac{1}{6}\beta_1^3\beta_2\nonumber\\
&=&
\displaystyle\frac{(b-c)^5-(a-c)^5}{5!}-\displaystyle\frac{\beta_1^3\beta_2A_5}{6}+
\underbrace{\displaystyle\frac{(b-c)^6-(a-c)^6}{6\beta_1^5}}_{A_5}\cdot\displaystyle\frac{1}{6}\beta_1^3\beta_2\nonumber\\
&=&\displaystyle\frac{(b-c)^5-(a-c)^5}{5!}
\end{eqnarray}
and
\begin{eqnarray}\label{eq31}
&&\mbox{the coefficient of $f^{(5)}(c)$}\nonumber\\
&=&
\displaystyle\frac{(b-c)^6-(a-c)^6}{6\beta_1^5}\cdot
\displaystyle\frac{1}{5!}\beta_1^5=\displaystyle\frac{(b-c)^6-(a-c)^6}{6!}.
\end{eqnarray}

\medskip
Using  (\ref{eq27}), (\ref{eq28}),  (\ref{eq29}), (\ref{eq30}) and (\ref{eq31}), we get from  (\ref{eq26}) and 
 (\ref{eq24}) that
\begin{eqnarray*}\label{eq}
&&\Big(\Gamma_{\beta_1, \dots, \beta_5}\circ\Lambda(f)\circ\Omega_{\beta_1, \dots, \beta_5}\Big)(c)\\
&=&
(b-a) f(c)+\displaystyle\sum_{i=1}^5\displaystyle\frac{(b-c)^{i+1}-(a-c)^{i+1}}{(i+1)!}f^{(i)}(c)\\
&=&\displaystyle\sum_{i=0}^5\displaystyle\frac{(b-c)^{i+1}-(a-c)^{i+1}}{(i+1)!}f^{(i)}(c)
=\displaystyle\int_a^b T_{5; c}(x)dx
\end{eqnarray*}

\medskip
This completes the proof of Proposition \ref{pr2.1}.

\hfill\raisebox{1mm}{\framebox[2mm]{}}

\bigskip
\section{An Example} 

\medskip
For convenience, the $5$th-order automatic integration  $\big(\Lambda,\, \Omega_{\beta_1, \dots, \beta_5}, \Gamma_{\beta_1, \dots, \beta_5}\big)$ induced by $\mathcal{R}^{(6)}$  
 in Proposition \ref{pr2.1}  will be also denoted by the 
 {\bf   $\Big(\mathcal{R}^{(6)}_{ \beta_1, \dots, \beta_5},\, T_{5; c}(x) \Big)$-automatic integration}.  The different choices of the parameters $\beta_1\ne 0$, $\beta_2, \dots , \beta_5$ give different ways of doing automatic  integration to approximate $\displaystyle\int_a^b f(x) dx$  by computing 
 $\displaystyle\int_a^b T_{5; c}(x) dx$ exactly. 
 
 \bigskip
After denoting a  $\mathcal{R}^{(6)}$-number $x+\displaystyle\sum_{i=1}^5a_i\varepsilon^i
\in \mathcal{R}^{(6)}$  by a $6$-tuple $(x, a_1,  \dots ,  a_5)$ of real numbers,  the algorithm, which compute  $\displaystyle\int_a^b T_{5; c}(x) dx$ exactly, can be written in a pseudo code as follows:

\bigskip
$\underline{\overline{\mbox{{\bf Algorithm}\quad $\Big(\mathcal{R}^{(6)}_{ \beta_1, \dots, \beta_5},\, T_{5; c}(x) \Big)$-automatic integration}}}$

\begin{itemize}
\item {\bf Input:} Three real numbers $a$, $b$,  $c$ and an analytic function  
$f:\in  {\bf D}^{\infty}[a, \, b]$.
\item {\bf Output:} A  real number $z$.
\end{itemize}

\begin{enumerate}
\item Start.
\item Get the $\mathcal{R}^{(6)}$-extension $\overline{f}: \mathcal{R}^{(6)}\to \mathcal{R}^{(6)}$.
\item Compute the value of the function $\overline{f}$ at $(c, \beta_1, \dots, \beta_5)$ to get the $6$-tuple 
$(f(c), y_1, \dots , y_5)$ of real numbers.
\item Compute  $A_5:=\displaystyle\frac{(b-c)^6-(a-c)^6}{6\beta_1^5}$.
\item Compute  $A_4:=\displaystyle\frac{(b-c)^5-(a-c)^5}{5\beta_1^4}-\displaystyle\frac{4\beta_2A_5}{\beta_1}$.
\item Compute  $A_3:=\displaystyle\frac{(b-c)^4-(a-c)^4}{4\beta_1^3}-\displaystyle\frac{3\beta_2A_4}{\beta_1}
-3\left(\displaystyle\frac{\beta_3}{\beta_1}+\displaystyle\frac{\beta_2^2}{\beta_1^2}\right)A_5$.
\item Compute $A_2:=\displaystyle\frac{(b-c)^3-(a-c)^3}{3\beta_1^2}-\displaystyle\frac{2\beta_2A_3}{\beta_1}
-\left(\displaystyle\frac{2\beta_3}{\beta_1}+\displaystyle\frac{\beta_2^2}{\beta_1^2}\right)A_4+$

 \qquad\qquad\qquad $-2\left(\displaystyle\frac{\beta_4}{\beta_1}+\displaystyle\frac{\beta_2\beta_3}{\beta_1^2}\right)A_5$.
\item Compute $A_1:=\displaystyle\frac{(b-c)^2-(a-c)^2}{2\beta_1}-\displaystyle\frac{\beta_2A_2}{\beta_1}
-\displaystyle\frac{\beta_3A_3}{\beta_1}-\displaystyle\frac{\beta_4A_4}{\beta_1}-\displaystyle\frac{\beta_5A_5}{\beta_1}$.
\item Compute $z=(b-a)f(c)+\displaystyle\sum_{i=1}^5 A_i y_i.$
\item Display  $z$.
\item Stop.
\end{enumerate}

 \medskip
 In the remaining of this section, we will choose
 $\beta_1=1$, $\beta_2=\beta_3=\beta_4=\beta_5= 0$  and  approximate $\displaystyle\int_0^2 e^{x^2} dx$ using  $5$ different ways which are based on the 
 $\Big(\mathcal{R}^{(6)}_{ 1, 0, 0, 0, 0},\, T_{5; c}(x) \Big)$-automatic integration, where 
 $T_{5; c}(x)$ is the the  $5$th-order Taylor polynomial for $f(x)=e^{x^2}$ with its  center at $c$. .
 
 \medskip
 Using  $\mathcal{R}^{(6)}$-extensions of the exponential function $e^x$, we have
 \begin{eqnarray}\label{eq32}
 &&e^{-c}\cdot \overline{\exp} (c+a_1\varepsilon +a_2\varepsilon^2)=1+
a_1 \varepsilon +\left(\displaystyle\frac{a_1^2}{2}+a_2\right)\varepsilon^2+\left(\displaystyle\frac{a_1^3}{2}+a_1a_2\right)\varepsilon^3+\nonumber\\
&&\quad 
+\left(\displaystyle\frac{a_1^4}{24}+\displaystyle\frac12 a_1^2a_2
+\displaystyle\frac12 a_2^2\right)\varepsilon^4+\Big(\displaystyle\frac{1}{120}a_1^5 +\displaystyle\frac{1}{6}a_1^3a_2+\displaystyle\frac{a_1a_2^2}{2}\Big)\varepsilon^5,
 \end{eqnarray}
 where $c$, $a_1$ and $a_2\in \mathcal{R}$.
 
\bigskip
{\bf Way 1:  $5$-automatic approximation of  $\displaystyle\int_0^2 e^{x^2}  dx$ at its {\bf center}  $0$.}

\medskip
By (\ref{eq6}) and (\ref{eq32}), we have
\begin{equation}\label{eq33}
 e^{-0}\cdot \overline{\exp}\Big(\big (\Omega_{1, 0,0,0,0}(0)\big)^2\Big)= \overline{\exp}\Big((0+\varepsilon)^2\Big)=\overline{\exp}(\varepsilon^2)=1
 +\varepsilon^2+\displaystyle\frac12 \varepsilon^4.
\end{equation}

\medskip
Using (\ref{eq8}), we have
\begin{equation}\label{eq34}
 A_2=\displaystyle\frac{2^3-0^3}{3}=2.666666667 \qquad\mbox{and}\qquad A_4=\displaystyle\frac{2^5-0^5}{5}=6.4.
\end{equation}

\medskip
It follows from  (\ref{eq7}) , (\ref{eq33})  and (\ref{eq34}) that
\begin{eqnarray*}
&&\displaystyle\int_0^2 e^{x^2}  dx\approx\displaystyle\int_0^2 T_{5; 0}(x) dx
=\Big(\Gamma_{1, 0,0,0,0}\circ (\Lambda(e^{x^2}))\circ\Omega_{1, 0,0,0,0}\Big)(0)\nonumber\\
&=&\Gamma_{1, 0,0,0,0}\left( \overline{\exp}\Big(\big (\Omega_{1, 0,0,0,0}(0)\big)^2\Big)\right)
=\Gamma_{1, 0,0,0,0}\left(1+\varepsilon^2+\displaystyle\frac12 \varepsilon^4\right) \nonumber\\
&=&(2-0)\times 1+2.666666667 \times 1+6.4\times\displaystyle\frac12=
7.866666667,
\end{eqnarray*}
which underestimates $\displaystyle\int_0^2 e^{x^2}  dx$ with the error $8.585961133$.
 
\bigskip
{\bf Way 2:  $5$-automatic approximation of  $\displaystyle\int_0^2 e^{x^2}  dx$ at its {\bf center}  $0.9$.}

\medskip
By (\ref{eq6}) and (\ref{eq32}), we have
\begin{eqnarray}\label{eq36}
 &&e^{-0.81}\cdot \overline{\exp}\Big( \big(\Omega_{1, 0,0,0,0}(0.9)\big)^2\Big)=e^{-0.81}\cdot  \overline{\exp}\Big((0.9+\varepsilon)^2\Big)\nonumber\\
 &=&e^{-0.81}\cdot \overline{\exp}(0.81+1.8\varepsilon+\varepsilon^2)
 =1+1.8\varepsilon+\left(\displaystyle\frac{1.8^2}{2}+1 \right)\varepsilon^2+\nonumber\\
 && +\left(\displaystyle\frac{1.8^3}{6}+1.8 \right)\varepsilon^3+\left(\displaystyle\frac{1.8^4}{24} +\displaystyle\frac{1.8^2}{2}+\displaystyle\frac12\right)\varepsilon^4+
\left(\displaystyle\frac{1.8^5}{120}+\displaystyle\frac{1.8^3}{6}+\displaystyle\frac{1.8}{2}\right)\varepsilon^5\nonumber\\
 &=&1+1.8\varepsilon+2.62\varepsilon^2+2.772\varepsilon^3+2.5574\varepsilon^4+2.029464\varepsilon^5.
\end{eqnarray}

\medskip
Using $b-c=2-0.9=1.1$,  $a-c=0-0.9=-0.9$ and (\ref{eq8}), we have
\begin{equation}\label{eq37}
\left\{\begin{array}{l}
 A_5=\displaystyle\frac{1.1^6-(-0.9)^6}{6}=\displaystyle\frac{1.771561-0.531441}{6}=0.206686666, \\\\
A_4=\displaystyle\frac{1.1^5-(-0.9)^5}{5}=\displaystyle\frac{1.61051+0.59049}{5}=0.4402, \\\\
A_3=\displaystyle\frac{1.1^4-(-0.9)^4}{4}=\displaystyle\frac{1.4641-0.6561}{4}=0.202, \\\\
A_2=\displaystyle\frac{1.1^3-(-0.9)^3}{3}=\displaystyle\frac{1.331+0.729}{3}=0.686666666, \\\\
A_1=\displaystyle\frac{1.1^2-(-0.9)^2}{2}=\displaystyle\frac{1.21-0.81}{2}=0.2.
\end{array}\right.
\end{equation}

\medskip
It follows from  (\ref{eq7}) , (\ref{eq36})  and (\ref{eq37}) that
\begin{eqnarray*}
&&\displaystyle\int_0^2 e^{x^2}  dx\approx\displaystyle\int_0^2 T_{5;\, 0.9}(x) dx
=e^{0.81}\Big\{(2-0)\times 1+0.2\times 1.8+\nonumber\\
&&\qquad +0.686666666\times 2.62+0.202\times 2.772+ 0.4402\times 2.5574+\nonumber\\
&& +0.206686666\times2.029464 \Big\}=e^{0.81}\Big\{2+0.36+1.799066665+0.559944+\nonumber\\
&&\qquad +1.12576748+0.419463147\Big\}=e^{0.81}\times 6.264241293=14.081438,
\end{eqnarray*}
which underestimates $\displaystyle\int_0^2 e^{x^2}  dx$ with the error $2.371189803$.

\medskip
{\bf Remark}  Although both Way 1 and Way 2 use one subinterval to evaluate the $5$-automatic approximations of  the definite integral $\displaystyle\int_0^2 e^{x^2} dx$,  the accuracy of  the $5$-automatic approximation in Way 2 is much better than the accuracy of  the $5$-automatic approximation in Way 1 after the center is changed from $0$ to $0.9$.  Hence, the accuracy of the automatic approximations of  a definite integral  not only depends on the the number of the subintervals, but also depends on the choices of the centers.

\bigskip
{\bf Way 3:  ${5, 5\choose 1.38}$-automatic approximation of  $\displaystyle\int_0^2 e^{x^2}  dx$ at its {\bf center}  

\qquad\quad\, $(0,\, 1.38)$.}

\medskip
First , we compute $5$-automatic approximation of  $\displaystyle\int_0^{1.38} e^{x^2}  dx$ at its  center $0$.

\medskip
Using (\ref{eq8}), we have
\begin{equation}\label{eq38}
 A_2=\displaystyle\frac{1.38^3-0^3}{3}=0.876024 \quad\mbox{and}\quad A_4=\displaystyle\frac{1.38^5-0^5}{5}=1.000980063.
\end{equation}

\medskip
It follows from  (\ref{eq7}) , (\ref{eq33})  and (\ref{eq38}) that
\begin{eqnarray}\label{eq39}
&&\displaystyle\int_0^{1.38} e^{x^2}  dx\approx\displaystyle\int_0^{1.38} T_{5; 0}(x) dx 
=(1.38-0)\times 1+ \nonumber\\ 
&&\qquad +0.876024\times 1+1.000980063\times\displaystyle\frac12=2.756514032.
\end{eqnarray}

\medskip
Next , we compute $5$-automatic approximation of  $\displaystyle\int_{1.38}^2 e^{x^2}  dx$ at its  center $1.38$.

\medskip
By (\ref{eq6}) and (\ref{eq32}), we have
\begin{eqnarray}\label{eq40}
 &&e^{-1.38^2}\cdot \overline{\exp}\Big( \big(\Omega_{1, 0,0,0,0}(1.38)\big)^2\Big)=e^{-1.38^2}\cdot  \overline{\exp}\Big((1.38+\varepsilon)^2\Big)\nonumber\\
 &=&e^{-1.9044}\cdot \overline{\exp}(1.9044+2.76\varepsilon+\varepsilon^2)
 =1+2.76\varepsilon+\left(\displaystyle\frac{2.76^2}{2}+1 \right)\varepsilon^2+\nonumber\\
 && +\left(\displaystyle\frac{2.76^3}{6}+2.76 \right)\varepsilon^3+\left(\displaystyle\frac{2.76^4}{24} +\displaystyle\frac{2.76^2}{2}+\displaystyle\frac12\right)\varepsilon^4+\nonumber\\
&&\qquad +\left(\displaystyle\frac{2.76^5}{120}+\displaystyle\frac{2.76^3}{6}+\displaystyle\frac{2.76}{2}\right)\varepsilon^5
 =1+2.76\varepsilon+\nonumber\\
 &&\quad +4.8088\varepsilon^2+6.264096\varepsilon^3+6.72662624\varepsilon^4+6.218736084\varepsilon^5.
\end{eqnarray}

\medskip
Using $b-c=2-1.38=0.62$,  $a-c=1.38-1.38=0$ and (\ref{eq8}), we have
\begin{equation}\label{eq41}
\left\{\begin{array}{l}
 A_5=\displaystyle\frac{0.62^6-0^6}{6}=\displaystyle\frac{0.056800235}{6}=0.009466705, \\\\
A_4=\displaystyle\frac{0.62^5-0^6}{5}=\displaystyle\frac{0.116029062}{5}=0.018322656, \\\\
A_3=\displaystyle\frac{0.62^4-0^6}{4}=\displaystyle\frac{0.14776336}{4}=0.03694084, \\\\
A_2=\displaystyle\frac{0.62^3-0^6}{3}=\displaystyle\frac{0.238328}{3}=0.079442666, \\\\
A_1=\displaystyle\frac{0.62^2-0^6}{2}=\displaystyle\frac{0.3844}{2}=0.1922. 
\end{array}\right.
\end{equation}

\medskip
It follows from  (\ref{eq7}) , (\ref{eq40})  and (\ref{eq41}) that
\begin{eqnarray}\label{eq42}
&&\displaystyle\int_{1.38}^2 e^{x^2}  dx\approx\displaystyle\int_{1.38}^2 T_{5;\, 1.38}(x) dx
=e^{1.9044}\Big\{(2-1.38)\times 1+0.1922\times 2.76+\nonumber\\
&&\qquad +0.079442666\times 4.8088 +0.03694084\times 6.264096+\nonumber\\
&&\quad + 0.018322656\times 6.72662624+0.009466705\times  6.218736084\Big\}\nonumber\\
&=&e^{1.9044}\Big\{0.62+0.530472 +0.382023892+0.231400968+0.123249658+\nonumber\\
&&\qquad +0.0588570939\Big\}=e^{1.9044}\times 1.946017458 =13.06824125.
\end{eqnarray}

\medskip
By  (\ref{eq39})  and (\ref{eq42}), we get
\begin{eqnarray*}
&&\displaystyle\int_0^2 e^{x^2}  dx
=\displaystyle\int_0^{1.38} e^{x^2}  dx+\displaystyle\int_{1.38}^2 e^{x^2}  dx\\
&\approx&
\displaystyle\int_0^{1.38} T_{5; 0}(x) dx +\displaystyle\int_{1.38}^2 T_{5;\, 1.38}(x) dx\\
&=&2.756514032+13.06824125=15.82475528,
\end{eqnarray*}
which underestimates $\displaystyle\int_0^2 e^{x^2}  dx$ with the error $0.627872518$.

\medskip
Let $T(n)$ be the Trapezoid Rule approximation to $\displaystyle\int_0^2 e^{x^2}  dx$ using $n$ subintervals. Then we have
$
\displaystyle\int_0^2 e^{x^2}  dx\approx T(8)=17.5650858,
$
which overestimates $\displaystyle\int_0^2 e^{x^2}  dx$ with the error $1.112458$. Hence,   the 
${5, 5\choose 1.38}$-automatic approximation  to $\displaystyle\int_0^2 e^{x^2}  dx$ at its  center 
$(0, 1.38)$, which uses two subintervals,  is better than 
the Trapezoid Rule approximation to $\displaystyle\int_0^2 e^{x^2}  dx$ using $8$ subintervals.

\bigskip
{\bf Way 4:  ${5, 5\choose 1.38}$-automatic approximation of  $\displaystyle\int_0^2 e^{x^2}  dx$ at its {\bf center}  

\qquad\quad\, $(0.65,\, 1.38)$.}

\medskip
Let us  compute $5$-automatic approximation of  $\displaystyle\int_0^{1.38} e^{x^2}  dx$ at its  center $0.65$. 

\medskip
By (\ref{eq6}) and (\ref{eq32}), we have
\begin{eqnarray}\label{eq43}
 &&e^{-0.65^2}\cdot \overline{\exp}\Big( \big(\Omega_{1, 0,0,0,0}(0.65)\big)^2\Big)=e^{-0.65^2}\cdot  \overline{\exp}\Big((0.65+\varepsilon)^2\Big)\nonumber\\
 &=&e^{-0.4225}\cdot \overline{\exp}(0.4225+1.3\varepsilon+\varepsilon^2)
 =1+1.3\varepsilon+\left(\displaystyle\frac{1.3^2}{2}+1 \right)\varepsilon^2+\nonumber\\
 && +\left(\displaystyle\frac{1.3^3}{6}+1.3 \right)\varepsilon^3+\left(\displaystyle\frac{1.3^4}{24} +\displaystyle\frac{1.3^2}{2}+\displaystyle\frac12\right)\varepsilon^4+\nonumber\\
&&\qquad +\left(\displaystyle\frac{1.3^5}{120}+\displaystyle\frac{1.3^3}{6}+\displaystyle\frac{1.3}{2}\right)\varepsilon^5
 =1+1.3\varepsilon+1.845\varepsilon^2+\nonumber\\
 &&\quad +1.6661666667\varepsilon^3+1.464004167\varepsilon^4+1.047107749\varepsilon^5.
\end{eqnarray}

\medskip
Using $b-c=1.38-0.65=0.73$,  $a-c=0-0.65=-0.65$ and (\ref{eq8}), we have
$$
\left\{\begin{array}{l}
 A_5=\displaystyle\frac{0.73^6-(-0.65)^6}{6}=\displaystyle\frac{0.151334226-0.07541889}{6}=0.012652555, \\\\
 A_4=\displaystyle\frac{0.73^5-(-0.65)^5}{5}=\displaystyle\frac{0.207307159+0.116029062}{5}=0.064667244, \\\\
 A_3=\displaystyle\frac{0.73^4-(-0.65)^4}{4}=\displaystyle\frac{0.28398241-0.17850625}{4}=0.02636904, \\\\
 A_2=\displaystyle\frac{0.73^3-(-0.65)^3}{3}=\displaystyle\frac{0.389017+0.274625}{3}=0.221214, \\\\
 A_1=\displaystyle\frac{0.73^2-(-0.65)^2}{2}=\displaystyle\frac{0.5329-0.4225}{2}=0.0552.
\end{array}\right.
$$

By the values of $A_i$ with $1\le i\le 5$ above and (\ref{eq43}) , we get
\begin{eqnarray}\label{eq44}
&&\displaystyle\int_{0}^{1.38} e^{x^2}  dx\approx\displaystyle\int_{0}^{1.38} T_{5;\, 0.65}(x) dx
=e^{0.4225}\Big\{(1.38-0)\times 1+\nonumber\\
&&\qquad +0.0522\times 1.3+0.221214\times 1.845 +0.02636904\times 1.6661666667+\nonumber\\
&&\quad + 0.064667244\times 1.464004167+0.012652555\times 1.047107749 \Big\}\nonumber\\
&=&e^{0.4225}\Big\{1.38+0.07176+0.40813983+0.043935215+0.094673114+\nonumber\\
&&\qquad +0.013248588\Big\}=e^{0.4225}\times  2.011756747=3.069480545.
\end{eqnarray}

\medskip
By  (\ref{eq42})  and (\ref{eq44}), we get
\begin{eqnarray*}
&&\displaystyle\int_0^2 e^{x^2}  dx
=\displaystyle\int_0^{1.38} e^{x^2}  dx+\displaystyle\int_{1.38}^2 e^{x^2}  dx\\
&\approx&
\displaystyle\int_0^{1.38} T_{5; 0.65}(x) dx +\displaystyle\int_{1.38}^2 T_{5;\, 1.38}(x) dx\\
&=&3.069480545+13.06824125=16.13772199,
\end{eqnarray*}
which underestimates $\displaystyle\int_0^2 e^{x^2}  dx$ with the error $0.314906076$.

\medskip
Let $M(n)$ be the Midpoint Rule approximation to $\displaystyle\int_0^2 e^{x^2}  dx$ using $n$ subintervals. Then we have
$
\displaystyle\int_0^2 e^{x^2}  dx\approx M(8)=15.9056767,
$
which underestimates $\displaystyle\int_0^2 e^{x^2}  dx$ with the error $0.5469511$. Hence,   the 
${5, 5\choose 1.38}$-automatic approximation  to $\displaystyle\int_0^2 e^{x^2}  dx$ at its  center 
 $(0.65, \,1.38)$,
, which uses two subintervals,  is better than 
the Midpoint Rule approximation to $\displaystyle\int_0^2 e^{x^2}  dx$ using $8$ subintervals.

\bigskip
{\bf Way 5:  ${5,\, 5,\,  5\choose 1.38,\, 1.39}$-automatic approximation of  $\displaystyle\int_0^2 e^{x^2}  dx$ at its {\bf center}  

\qquad\quad\, $(0.65, \, 1.38,\, 1.69)$.}

\medskip
First, we  compute $5$-automatic approximation of  $\displaystyle\int_{1.38}^{1.39} e^{x^2}  dx$ at its  center $1.38$.

\medskip
Using $b-c=1.39-1.38=0.01$,  $a-c=1.38-1.38=0$ and (\ref{eq8}), we have
\begin{equation}\label{eq45}
\left\{\begin{array}{l}
 A_5=\displaystyle\frac{0.01^6-0^6}{6}=\displaystyle\frac{1}{6\times 10^{12}}, \\\\
A_4=\displaystyle\frac{0.01^5-0^6}{5}=\displaystyle\frac{1}{5\times 10^{10}}, \\\\
A_3=\displaystyle\frac{0.01^4-0^6}{4}=0.000000002, \\\\
A_2=\displaystyle\frac{0.01^3-0^6}{3}=0.000000333, \\\\
A_1=\displaystyle\frac{0.01^2-0^6}{2}=0.00005. 
\end{array}\right.
\end{equation}

\medskip
It follows from  (\ref{eq7}) , (\ref{eq40})  and (\ref{eq45}) that
\begin{eqnarray}\label{eq46}
&&\displaystyle\int_{1.38}^{1.39} e^{x^2}  dx\approx\displaystyle\int_{1.38}^{1.39} T_{5;\, 1.38}(x) dx
=e^{1.9044}\Big\{(1.39-1.38)\times 1+0.00005\times 2.76+\nonumber\\
&&\qquad +0.000000333\times 4.8088 +0.000000002\times 6.264096+\nonumber\\
&&\quad + \displaystyle\frac{1}{5\times 10^{10}}\times 6.72662624+\displaystyle\frac{1}{6\times 10^{12}}\times  6.218736084\Big\}\nonumber\\
&\approx&e^{1.9044}\Big\{0.01+0.000138 +0.000001601+0.000000012\Big\}\nonumber\\
&=&e^{1.9044}\times 
0.010139613 =0.068091329.
\end{eqnarray}

\medskip
Next , we compute $5$-automatic approximation of  $\displaystyle\int_{1.39}^2 e^{x^2}  dx$ at its  center $1.69$.

\medskip
By (\ref{eq6}) and (\ref{eq32}), we have
\begin{eqnarray}\label{eq47}
 &&e^{-1.69^2}\cdot \overline{\exp}\Big( \big(\Omega_{1, 0,0,0,0}(1.69)\big)^2\Big)=e^{-1.69^2}\cdot  \overline{\exp}\Big((1.69+\varepsilon)^2\Big)\nonumber\\
 &=&e^{-2.8561}\cdot \overline{\exp}(2.8561+3.38\varepsilon+\varepsilon^2)
 =1+3.38\varepsilon+\left(\displaystyle\frac{3.38^2}{2}+1 \right)\varepsilon^2+\nonumber\\
 && +\left(\displaystyle\frac{3.38^3}{6}+3.38 \right)\varepsilon^3+\left(\displaystyle\frac{3.38^4}{24} +\displaystyle\frac{3.38^2}{2}+\displaystyle\frac12\right)\varepsilon^4+\nonumber\\
&&\qquad +\left(\displaystyle\frac{3.38^5}{120}+\displaystyle\frac{3.38^3}{6}+\displaystyle\frac{3.38}{2}\right)\varepsilon^5
 =1+3.38\varepsilon+\nonumber\\
 &&\quad +6.7122\varepsilon^2+9.815745333\varepsilon^3+11.65040481\varepsilon^4+11.80197178\varepsilon^5.
\end{eqnarray}

\medskip
Using $b-c=2-1.69=0.31$,  $a-c=1.39-1.69=-0.3$ and (\ref{eq8}), we have
\begin{equation}\label{eq48}
\left\{\begin{array}{l}
 A_5=\displaystyle\frac{0.31^6-(-0.3)^6}{6}=\displaystyle\frac{0.000887503-0.000729}{6}
 =0.000026417, \\\\
 A_4=\displaystyle\frac{0.31^5-(-0.3)^5}{5}=\displaystyle\frac{0.002862915+0.00243}{5}
 =0.001058583, \\\\
 A_3=\displaystyle\frac{0.31^4-(-0.3)^4}{4}=\displaystyle\frac{0.00923521-0.0081}{4}
 =0.000283802, \\\\
 A_2=\displaystyle\frac{0.31^3-(-0.3)^3}{3}=\displaystyle\frac{0.029791+0.027}{3}
 =0.018930333, \\\\
 A_1=\displaystyle\frac{0.31^2-(-0.3)^2}{3}=\displaystyle\frac{0.0961-0.09}{2}
 =0.00305.
\end{array}\right.
\end{equation}

\medskip
It follows from  (\ref{eq7}) , (\ref{eq47})  and (\ref{eq48}) that
\begin{eqnarray}\label{eq49}
&&\displaystyle\int_{1.39}^2 e^{x^2}  dx\approx\displaystyle\int_{1.39}^2 T_{5;\, 1.69}(x) dx
=e^{2.8561}\Big\{(2-1.39)\times 1+0.00305\times 3.38+\nonumber\\
&&\qquad +0.018930333\times  6.7122+0.000283802\times 9.815745333+\nonumber\\
&&\quad +0.001058583 \times 11.65040481+0.000026417\times  11.80197178\Big\}\nonumber\\
&=&e^{2.8561}\Big\{0.61+0.010309+0.127064181+0.002785728+0.01233292+\nonumber\\
&&\qquad +0.000311772\Big\}=e^{2.8561}\times  0.762803601=13.26786991.
\end{eqnarray}

\medskip
By  (\ref{eq44}),  (\ref{eq46})  and (\ref{eq49}), we get
\begin{eqnarray}\label{eq50}
&&\displaystyle\int_0^2 e^{x^2}  dx
=\displaystyle\int_0^{1.38} e^{x^2}  dx+\displaystyle\int_{1.38}^{1.39} e^{x^2}  dx+
\displaystyle\int_{1.39}^{2} e^{x^2}  dx\nonumber\\
&\approx&
\displaystyle\int_0^{1.38} T_{5; 0.65}(x) dx +\displaystyle\int_{1.38}^{1.39} T_{5;\, 1.38}(x) dx+\displaystyle\int_{1.39}^2 T_{5;\, 1.69}(x) dx\nonumber\\
&=&3.069480545+0.068091329+13.26786991=16.40544197,
\end{eqnarray}
which underestimates $\displaystyle\int_0^2 e^{x^2}  dx$ with the error $0.047186016$.

\medskip
Let $S(n)$ be Simpson Rule approximation to $\displaystyle\int_0^2 e^{x^2}  dx$ using $n$ subintervals. Then we have
$
\displaystyle\int_0^2 e^{x^2}  dx\approx S(8)=16.5385947,
$
which overestimates $\displaystyle\int_0^2 e^{x^2}  dx$ with the error $0.0859669$. Hence,  
the ${5,\, 5,\,  5\choose 1.38,\, 1.39}$-automatic approximation of  $\displaystyle\int_0^2 e^{x^2}  dx$ at its center $(0.65, \, 1.38,\, 1.69)$,  which uses $3$ subintervals,  is better than 
Simpson Rule approximation to $\displaystyle\int_0^2 e^{x^2}  dx$ using $8$ subintervals.

\bigskip
By the way, one can check that
$$
\displaystyle\int_0^2 e^{x^2}  dx\approx M(16)=16.3118539 
$$
which underestimates $\displaystyle\int_0^2 e^{x^2}  dx$ with the error $0.1407739$ and
$$
\displaystyle\int_0^2 e^{x^2}  dx\approx T(16)=16.7353812
$$
which overestimates $\displaystyle\int_0^2 e^{x^2}  dx$ with the error $0.2827535$. Hence,  
the ${5,\, 5,\,  5\choose 1.38,\, 1.39}$-automatic approximation of  $\displaystyle\int_0^2 e^{x^2}  dx$ at its center $(0.65, \, 1.38,\, 1.69)$,  which is given by  (\ref{eq50}) and  just uses $3$ subintervals,  is better than 
both the Midpoint Rule and the Trapezoid Rule approximation to $\displaystyle\int_0^2 e^{x^2}  dx$ using $16$ subintervals.

\end{document}